\documentclass{amsart}
\usepackage{amsmath}
\usepackage{graphicx}
\usepackage{amsfonts}
\usepackage{amssymb}

\newtheorem{theorem}{Theorem}
\theoremstyle{plain}

\newtheorem{definition}{Definition}
\newtheorem{example}{Example}

\newtheorem{proposition}{Proposition}
\newtheorem{remark}{Remark}

\numberwithin{equation}{section}

\begin{document}
\title{Filtered Hirsch algebras}
\author{Samson Saneblidze}
\address{A. Razmadze Mathematical Institute\\
Department of Geometry and Topology\\
M. Aleksidze st., 1\\
0193 Tbilisi, Georgia}
\email{sane@rmi.ge}
\thanks{This research described in this publication was made possible in
part by the grant \\
GNF/ST06/3-007 of the Georgian National Science Foundation}
\subjclass[2000]{Primary 55P35, 55S30; Secondary 55S05, 55U20}
\keywords{Hirsch algebra, filtered model, multiplicative resolution,
symmetric Massey product, Steenrod operation, Hochschild cohomology}
\date{}

\begin{abstract}
Motivated by the cohomology theory of loop spaces, we consider a special
class of higher order homotopy commutative differential graded algebras and
construct the filtered Hirsch model for such an algebra $A$. When $x\in H(A)$
with $\mathbb{Z}$ coefficients and $x^{2}=0,$ the symmetric Massey products $%
\langle x\rangle ^{n}$ with $n\geq 3$ have a finite order (whenever
defined). \ However, if $\Bbbk $ is a field of characteristic zero, $\langle
x\rangle ^{n}$ is defined and vanishes in $H(A\otimes \Bbbk )$ for all $n$.
If $p$ is an odd prime, the Kraines formula $\langle x\rangle ^{p}=-\beta
\mathcal{P}_{1}(x)$ lifts to $H^{\ast }(A\otimes {\mathbb{Z}}_{p}).$
Applications of the existence of polynomial generators in the loop homology
and the Hochschild cohomology with a $G$-algebra structure are given.
\end{abstract}

\maketitle
\tableofcontents

\section{Introduction}

In this paper we investigate a special class of homotopy commutative
algebras called \emph{Hirsch algebras} \cite{KSpermu}. When the structural
operations of a Hirsch algebra $A$ agree component-wise with those of a
homotopy $G$-algebra (HGA), the pre-Jacobi axiom can fail \cite{Voronov},
\cite{GJ}, \cite{KScubi}, \cite{Voronov2} and the induced product on the bar
construction $BA$ is not necessarily associative. Indeed, the theory of loop
space cohomology suggests that it is impossible in general, to construct a
small model for $H^{\ast }\left( \Omega X\right) $ in the category of HGAs.
The investigation here applies a perturbation theory that extends the
well-developed perturbation theories for differential graded modules and
differential graded algebras (dgas) \cite{berika3}, \cite{Gugenheim}, \cite%
{hueb-kade}, \cite{hal-sta}, \cite{sane}, \cite{saneDerived}.

One difficulty encountered when constructing a theory of homological algebra
for Hirsch algebras is that the Steenrod cochain product $a\smile _{1}b$
fails to be a cocycle even for cocycles $a$ and $b.$ Consequently $a\smile
_{1}b$ does not necessarily lift to cohomology. We control such difficulties
by introducing the notion of a \emph{filtered} Hirsch algebra, which can be
thought of as a specialization of a distinguished resolution in the sense of
\cite{gug-may} (see also \cite{HMS}). On the other hand, the filtered Hirsch
model $(RH,d+h)$ of a Hirsch algebra $A$ is itself a Hirsch algebra whose
structural operations $E_{p,q}:RH^{\otimes p}\otimes RH^{\otimes
q}\longrightarrow RH$ are completely determined by the commutative graded
algebra (cga) structure of $H=H(A,d_{A});$ furthermore, the perturbation $
h:RH\rightarrow RH$ of the resolution differential $d$ is determined by the
Hirsch algebra structure on $A$ (Theorem \ref{filtered}). Thus by ignoring
the operations $E_{p,q}$ we obtain a multiplicative resolution $
(RH,d)\rightarrow (H,0)$ of the cga $H$ thought of as a non-commutative
version of its Tate-Jozefiak resolution (\cite{Tate}, \cite{Jozefiak}) and
the filtered model of the dga $A$ is the perturbation $(RH,d+h)\rightarrow
(A,d_{A})$ in \cite{sane} (such a filtered model in the category of cdgas
over a field of characteristic zero was constructed by Halperin and Stasheff
in \cite{hal-sta}).

A Hirsch resolution always admits a binary operation $\cup _{2},$ which can
be viewed as \emph{divided }Steenrod $\smile _{2}$-operation. This leads to
the notion of a \emph{quasi-homotopy commutative} Hirsch algebra (QHHA)
introduced here. We note that in general, the construction of a Hirsch map $
(RH,d+h)\rightarrow A$ compatible with a QHHA structure on $A$ is obstructed
by the non-free action of $Sq_{1}$ on its cohomology $H(A).$

Every cdga $H$ can be thought of as a trivial Hirsch algebra in which the
operations $E_{p,q}\equiv 0$ for all $p,q\geq 1$. However, we exhibit an
example of a cohomology algebra $H=H(A)$ with a non-trivial Hirsch algebra
structure determined by $Sq_{1}.$

For a Hirsch algebra $A$ over the integers, we establish some formulas
relating the structural operations $E_{p,q}$ with syzygies in $(RH,d)$ that
arise from a single element $x\in H(A)$ with $x^{2}=0.$ Whereas the $n$-fold
symmetric Massey product $\langle x\rangle ^{n}$ with $n\geq 3$ is defined
in $H(A)$ (\cite{kraines}, \cite{kochman}), our formulas imply that $\langle
x\rangle ^{n}$ has finite order. Note that when $A$ is an algebra over a
field $\Bbbk $ of characteristic zero, $\langle x\rangle ^{n}$ is defined
and vanishes for all $n\geq 3$ (Theorem \ref{mzero}). As a consequence we
have (compare \cite{browder}):

\noindent \textbf{Theorem A.} \emph{\ Let }$X$\emph{\ be a simply connected
space, let }$\Bbbk $\emph{\ be a field of characteristic zero and let }$
\sigma _{\ast }:H_{\ast }(\Omega X;\Bbbk )\rightarrow H_{\ast +1}(X;\Bbbk )$
\emph{be the suspension map. If }$y\notin \operatorname{Ker}\sigma _{\ast }$ \emph{
and }$y^{2}\neq 0$\emph{, then }$y^{n}\neq 0$\emph{\ for all }$n\geq 2.$

Given an odd prime $p,$ consider the Hirsch algebra $A\otimes {\mathbb{Z}}
_{p},$ let $x\in H^{2m+1}(A\otimes {\mathbb{Z}}_{p}),$ and let $\beta $ be
the Bockstein operator. We obtain the formula
\begin{equation}
\langle x\rangle ^{p}=-\beta \mathcal{P}_{1}(x),  \label{one}
\end{equation}
which has the same form as Kraines's formula in \cite{kraines}, however, the
cohomology operation $\mathcal{P}_{1}:H^{2m+1}(A\otimes {\mathbb{Z}}
_{p})\rightarrow H^{2mp+1}(A\otimes {\mathbb{Z}}_{p})$ in (\ref{one}) is
canonically determined by the iteration of the $\smile _{1}$-product on $
A\otimes {\mathbb{Z}}_{p}$ (Theorem \ref{krainest}). Dually, if $A$ is the
singular chains on the triple loop space $\Omega ^{3}X$, we can identify $
\mathcal{P}_{1}$ with the Dyer-Lashof operation (see \cite{kochman}).
In fact the validity of (\ref{one}) in a general algebraic framework  is conjectured by May
\cite[Section 6]{may}.
Furthermore, when $X=BF_{4},$ the classifying space of the exceptional group $F_{4},$ we
exhibit explicit perturbations in the filtered model of $X$ and recover
formula (\ref{one}) in $H^{\ast }(X;{\mathbb{Z}}_{3})$.

Although Theorem 1 provides a theoretical model of a Hirsch algebra $A$
endowed with higher order operations $E_{p,q},$ in practice one can
construct a small
\emph{multiplicative}
 model for recognizing $H^{\ast }(BA)$
as an algebra in which the product is determined only by the binary
operation $E_{1,1}=\,\smile _{1}$. Thus, a (minimal) multiplicative
resolution of $H^{\ast }(A)$ endowed with a $\smile _{1}$-product provides
an economical way to calculate the algebra $H^{\ast }(BA).$
We apply this technique to
 the Hochschild cochain complex $A=C^{\bullet }(P;P)$
of an associative algebra $P$ over a field $\Bbbk $ of characteristic zero to establish the following

\noindent \textbf{Theorem B.} \emph{If the Hochschild cohomology }$H^{\ast
}=H(C^{\bullet }(P;P))$\emph{\ is a free algebra, then the Lie algebra
structure on }$Tor_{\ast }^{A}(\Bbbk ,\Bbbk )$\emph{\ is completely
determined by that of the }$G$\emph{-algebra }$H^{\ast }$\emph{$.$ Consequently, the
product }$\mu^{\ast} $\emph{\ on} $Tor_{\ast }^{A}(\Bbbk ,\Bbbk )$
\emph{is commutative if and only if the }$G$\emph{-product on }$H^{\ast }$\emph{is trivial.}

Some applications of filtered Hirsch algebras considered in an earlier
version of this paper are also considered in \cite{sanePOL}, \cite{saneBetti}
(see also \cite{saneFREE}, \cite{saneGeodesic}).

I wish to thank Jim Stasheff for helpful comments and suggestions.  I am also
indebted to the referee for a number of helpful comments and for having
suggested many improvements of the exposition.

\section{The category of Hirsch algebras}

\label{Hiralgebras} This section defines the generalized notion of a Hirsch
algebra applied here, the morphisms between them, and the notion of a Hirsch
resolution.

Let $\Bbbk $ be a commutative ring with unity $1$ and characteristic $\nu\, ;$
in the applications, $\Bbbk $ will be the integers $\mathbb{Z},$ a finite
field $\mathbb{Z}_{p}=\mathbb{Z}/p\mathbb{Z}$ with $p$ prime, or a field of
characteristic zero. Graded $\Bbbk $-modules $A^{\ast }$ are assumed to be
graded over ${\mathbb{Z}}.$ A module $A^{\ast }$ is connected if $%
A^{0}=\Bbbk ,$ and a non-negatively graded, connected module $A^{\ast }$ is $%
1$-\emph{reduced} if $A^{1}=0.$

For a module $A,$ let $T(A)=\bigoplus_{i=0}^{\infty }A^{\otimes i},$ where $%
A^{0}=\Bbbk ,$ be the tensor module of $A$. An element $a_{1}\otimes \cdots
\otimes a_{n}\in A^{\otimes n}$ is denoted by $[a_{1}|\cdots |a_{n}]$ when $%
T(A)$ is viewed as the tensor coalgebra or by $a_{1}\cdots a_{n}$ when $T(A)$
is viewed as the tensor algebra. We denote by $s^{-1}A$ the desuspension of $%
A$, i.e., $(s^{-1}A)^{i}=A^{i+1}$.

A dga $(A,d_{A})$ is assumed to be supplemented; in particular, it has the
form $A=\tilde{A}\oplus \Bbbk .$ The (reduced) bar construction $BA$ on $A$
is the tensor coalgebra $T(\bar{A}),\ \bar{A}=s^{-1}\tilde{A},$ with
differential $d=d_{1}+d_{2}$ given for $[\bar{a}_{1}|\dotsb |\bar{a}_{n}]\in
T^{n}(\bar{A})$ by
\begin{equation*}
d_{1}[\bar{a}_{1}|\dotsb |\bar{a}_{n}]=-\sum_{1\leq i\leq n}(-1)^{\epsilon
_{i-1}^{a}}[\bar{a}_{1}|\dotsb |\overline{d_{A}(a_{i})}|\dotsb |\bar{a}_{n}]
\end{equation*}%
and
\begin{equation*}
d_{2}[\bar{a}_{1}|\dotsb |\bar{a}_{n}]=-\sum_{1\leq i<n}(-1)^{\epsilon
_{i}^{a}}[\bar{a}_{1}|\dotsb |\overline{a_{i}a_{i+1}}|\dotsb |\bar{a}_{n}],
\end{equation*}%
where $\epsilon _{i}^{x}=|x_{1}|+\cdots +|{x_{i}}|+i.$

Let us generalize (slightly) the definition of a Hirsch algebra \cite%
{KSpermu}. Let $A$ be a dga and consider the dg module $({Hom}(BA\otimes
BA,A),\nabla )$, where $\nabla $ is the canonical $Hom$ differential. Since
the tensor product $BA\otimes BA$ is a dgc with the standard coalgebra
structure, the $\smile $-product induces a dga structure on $(Hom(BA\otimes
BA,A),\nabla ,\smile ).$

\begin{definition}
A \textbf{Hirsch algebra} is an associative dga $A$ equipped with
multilinear maps
\begin{equation*}
E_{p,q}:A^{\otimes p}\otimes A^{\otimes q}\rightarrow A,\ p,q\geq 0,\ p+q>0,
\end{equation*}
satisfying the following conditions:

\begin{enumerate}
\item[(i)] $\deg E_{p,q}=1-p-q$;

\item[(ii)] $E_{1,0}=Id=E_{0,1}\ \text{and}\ \ E_{p>1,0}=0=E_{0,q>1};$

\item[(iii)] The homomorphism $E:BA\otimes BA\rightarrow A$ defined by
\begin{equation}
E([\bar{a}_{1}|\dotsb |\bar{a}_{p}]\otimes \lbrack \bar{b}_{1}|\dotsb |\bar{b%
}_{q}])=E_{p,q}(a_{1},...,a_{p};b_{1},...,b_{q})  \label{two}
\end{equation}%
is a twisting cochain in the dga $(Hom(BA\otimes BA,A),\nabla ,\smile )$,
i.e., $\nabla E=-E\smile E.$
\end{enumerate}

A \emph{morphism} $f:A\rightarrow B$ between two Hirsch algebras is a dga
map $f$ that commutes with $E_{p,q}$ for all $p,q.$
\end{definition}

\noindent Condition (iii) implies that $\mu _{E}:BA\otimes BA\rightarrow BA$
is a chain map; thus $BA$ is a dg bialgebra whose multiplication $\mu _{E}
$ is not necessarily associative (compare \cite{GJ}, \cite{Voronov2}, \cite
{clark}, \cite{Khelaia}, \cite{munkholm}); in particular, $\mu _{_{{
E_{10}+E_{01}}}}$ is the shuffle product on $BA,$
and   a Hirsch algebra with $E_{p,q}\equiv 0$ for all $p,q\geq 1$ is
just a cdga (cf. (\ref{cup-one})).
 It is useful to express
equation (\ref{two}) component-wise:
\begin{multline}
dE_{p,q}(a_{1},...,a_{p};b_{1},...,b_{q})=\sum_{1\leq i\leq p}(-1)^{\epsilon
_{i-\!1}^{a}}E_{p,q}(a_{1},...,da_{i},...,a_{p};b_{1},...,b_{q})  \label{dif}
\\
\hspace{1.9in}+\sum_{1\leq j\leq q}^{{}}(-1)^{\epsilon _{p}^{a}+\epsilon
_{j-\!1}^{b}}E_{p,q}(a_{1},...,a_{p};b_{1},...,db_{j},...,b_{q}) \\
\hspace{1.9in}+\sum_{1\leq i<p}^{{}}\,(-1)^{\epsilon
_{i}^{a}}E_{p-1,q}(a_{1},...,a_{i}a_{i+1},...,a_{p};b_{1},...,b_{q}) \\
\hspace{1.87in}+\!\!\sum_{1\leq j<q}\!(-1)^{\epsilon _{p}^{a}+\epsilon
_{j}^{b}}E_{p,q-1}(a_{1},...,a_{p};b_{1},...,b_{j}b_{j+1},...,b_{q}) \\
\hspace{0.7in}+\!\!\sum_{\substack{ 0\leq i\leq p \\ 0\leq j\leq q \\ %
(i,j)\neq (0,0)}}\!\!\!\!\!(-1)^{\epsilon
_{i,j}}E_{i,j}(a_{1},\!...,a_{i};b_{1},\!...,b_{j})\cdot
E_{p-i,q-j}(a_{i+1},\!...,a_{p};b_{j+1},\!...,b_{q}), \\
{\epsilon _{i,j}}=\epsilon _{i}^{a}+\epsilon _{j}^{b}+(\epsilon
_{i}^{a}+\epsilon _{p}^{a})\epsilon _{j}^{b}+1.
\end{multline}
In particular, the operation $E_{1,1}$ satisfies conditions similar to
Steenrod's cochain $\smile _{1}$-product:
\begin{equation}\label{cup-one}
dE_{1,1}(a;b)-E_{1,1}(da;b)+(-1)^{|a|}E_{1,1}(a;db)=(-1)^{|a|}ab-(-1)^{|a|(|b|+1)}ba;
\end{equation}
consequently, $E_{1,1}$ measures the non-commutativity of the product $\cdot
$ on $A.$ We shall use the notation $a\smile _{1}b=E_{1,1}(a;b)$
interchangeably. The following special cases will also be important for us,
so we write them explicitly:

The Hirsch formulas up to homotopy
\begin{multline*}
dE_{2,1}(a,b;c)=E_{2,1}(da,b;c)-(-1)^{|a|}E_{2,1}(a,db;c)+(-1)^{|a|+|b|}E_{2,1}(a,b;dc)-
\\
(-1)^{|a|}(ab)\smile _{1}c+(-1)^{|a|+|b|+|b||c|}(a\smile
_{1}c)b+(-1)^{|a|}a(b\smile _{1}c)
\end{multline*}%
and
\begin{multline*}
dE_{1,2}(a;b,c)=E_{1,2}(da;b,c)-(-1)^{|a|}E_{1,2}(a;db,c)+(-1)^{|a|+|b|}E_{1,2}(a;b,dc)+
\\
(-1)^{|a|+|b|}a\smile _{1}(bc)-(-1)^{|a|+|b|}(a\smile
_{1}b)c-(-1)^{|a|(|b|-1)}b(a\smile _{1}c)
\end{multline*}%
tell us that the deviations of the binary operation $\smile _{1}$ from left
and right derivation of the $\cdot $ product are measured by the respective
boundaries of the operations $E_{1,2}$ and $E_{2,1}$ on three variables.

The following definition describes a class of Hirsch algebras  in which the $\smile_1$-product itself  is homotopy commutative (cf. (\ref{realcup2}) below).

\begin{definition}
\label{QHHA} A \textbf{quasi-homotopy commutative Hirsch algebra} (QHHA) is
a Hirsch algebra $A$ equipped with a binary product $\cup _{2}:A\otimes
A\rightarrow A$ such that
\begin{multline}
d(a\cup _{2}b)=da\cup _{2}b+(-1)^{|a|}a\cup _{2}db+(-1)^{|a|}a\smile
_{1}b+(-1)^{(|a|+1)|b|}b\smile _{1}a  \label{cup2} \\
-q(a;b),
\end{multline}
where $q(a;b)$ satisfies:

\begin{enumerate}
\item[\textit{$(2.4)_1$}] Leibniz rule: $dq(a;b)=-q(da;b)-(-1)^{|a|}q(a;db);$

\item[\textit{$(2.4)_2$}] Acyclicity: \ \, $[q(a,b)]=0\in H(A,d)$ for $
da=db=0.$
\end{enumerate}
\end{definition}

Note that $(2.4)_{1}$ follows from the equalities (2.2) and $d^{2}=0.$
Obviously, discarding the parameter $q(a;b),$ the above formula just becomes
the Steenrod formula for the $\smile _{2}$ -cochain product:
\begin{equation} \label{realcup2}
d(a\smile _{2}b)=da\smile _{2}b+(-1)^{|a|}a\smile _{2}db+(-1)^{|a|}a\smile
_{1}b+(-1)^{(|a|+1)|b|}b\smile _{1}a.
\end{equation}
However, $q(-;-)$ may be non-zero when passing to models constructed via
cohomology as below. In the following four  examples, the first is a
naturally occurring example of a \emph{cochain} Hirsch algebra (compare
Example \ref{hoh});   in the second example QHHA structures are considered for certain Hirsch algebras;   in the third and fourth
examples a Hirsch algebra
structure is lifted to the cohomology level. In fact, the fourth  example was
the original motivation for this paper.

\begin{example}
\label{top}The primary examples of Hirsch algebras for topological spaces $X$
are their cubical or simplicial cochain complexes \cite{KSpermu}, \cite%
{KScubi}, \cite{Khelaia}. In the simplicial case one can choose $E_{p,q}=0$
for $q\geq 2$ and obtain an HGA structure on the simplicial cochains $%
C^{\ast }(X;\Bbbk )$ \cite{baues3} (see also \cite{KScubi}). Furthermore,
the product $\mu _{E}$ on $BC^{\ast }(X;\Bbbk )$ gives the multiplicative
structure of the loop space cohomology $H^{\ast }(\Omega X;\Bbbk ).$

 Here
the cochain complex $C^{\ast }(X;\Bbbk )$ of a space $X$ is 1-reduced, since
by definition
 $C^{\ast }(X;\Bbbk )=C^{\ast }(\operatorname{Sing}^{1}X;\Bbbk )/C^{>0}(
\operatorname{Sing}\,x\, ;\Bbbk )$ where ${\operatorname{Sing}}^{1}X\subset {\operatorname{Sing}}X$ is
the Eilenberg 1-subcomplex generated by the singular simplices that send the
1-skeleton of the standard $n$-simplex $\Delta ^{n}$ to the base point $x$
of $X.$   Unlike the cubical cochains, the Hirsch algebra structure of the simplicial  cochains is \emph{associative}, i.e.,
the above product $\mu_{E}$ is associative.
\end{example}

\begin{example}\label{qhha}
First, note that the Hirsch algebras from the previous example  are also QHHA's by setting $\cup _{2}=\smile
_{2}$ and $q(-;-)=0.$
Let $A$  be a \emph{special Hirsch algebra}, i.e., $A$ is an associative Hirsch algebra and  $BA$ also admits a Hirsch algebra
structure.  Then $A$ is a QHHA since it
 admits a $\cup _{2}$-product satisfying (\ref{realcup2})
(cf. \cite{kade3}).
  An important example of a special Hirsch algebra
 is $A=C^{\ast }(X;\Bbbk )$ from the previous example (cf.
\cite{KSpermu}, \cite{SU2}).
Finally, for a QHHA  $A$  with $\nu$ to be zero or odd   and     $\smile_2$-product satisfying  (\ref{realcup2}), define the
\emph{divided} $\smile_2$-operation $\cup_2$ as
\[a\cup_2b =
\left\{
  \begin{array}{ll}
    \frac{1}{2}\, a\smile_2 a, & a=b   \\
    a\smile_2 b, & \text{otherwise}.
  \end{array}
\right.
\]
Then  $A$ with this  $\cup_2$-operation is again a QHHA.

\end{example}

\begin{example}
\label{freeH}
Let $(H,d=0)$ be a free cga  $H=S\langle\mathcal{H}^{\ast }\rangle$ generated by a graded set $\mathcal{H}^{\ast }.$ Then any map of sets $\tilde{E}_{p,q}:\mathcal{H}^{\times
p}\times \mathcal{H}^{\times q}\rightarrow H$ of degree $1-p-q$ extends to a
Hirsch algebra structure $E_{p,q}:H^{\otimes p}\otimes H^{\otimes
q}\rightarrow H$ on $H.$ Indeed, using formula (\ref{dif}) the construction
goes by induction on the sum $p+q.$ In particular, if only $\tilde{E}_{1,1}$
is non-zero then the image of $E_{p,q}$ for $p+q\geq 3$ is into the
submodule of $H$ spanned by the monomials of the form $\tilde{E}
_{1,1}(a_{1};b_{1})\cdots \tilde{E}_{1,1}(a_{k};b_{k})\cdot x$ for $
a_{i},b_{i}\in \mathcal{H},$ $x\in H,\,$and $k\geq 1.$
\end{example}

\begin{example}
\label{cartan}
 The argument in Example \ref{freeH} suggests how to lift a
Hirsch ${\mathbb{Z}}_{2}$-algebra structure from the cochain level to
cohomology. Given a Hirsch algebra $A,$ let $H=H^{\ast }(A).$ For a cocycle $%
a\in A^{m},$ one has $d_{A}E_{1,1}(a,a)=0$ and $Sq_{1}:H^{m}\rightarrow
H^{2m-1}$ is defined by
\begin{equation*}
\lbrack a]\rightarrow \lbrack E_{1,1}(a,a)].
\end{equation*}%
The trick here is to convert the Hirsch formulas up to homotopy on $A$ to
the Cartan formula $Sq_{1}(ab)=Sq_{1}a\!\cdot \!Sq_{0}b+Sq_{0}a\!\cdot
\!Sq_{1}b$ on $H\ $by fixing a set of multiplicative generators $\mathcal{H}%
\subset H.$ Define the map $\tilde{Sq}_{1,1}:\mathcal{H}\times \mathcal{H}%
\rightarrow H$ for $a,b\in \mathcal{H}$ by
\begin{equation*}
\tilde{Sq}_{1,1}(a;b)=
\begin{cases}
\,  Sq_{1}a, & a=b, \\
\,   0, & \text{otherwise}
\end{cases}
\end{equation*}
and extend to the operation $Sq_{1,1}:H\otimes H\rightarrow H$ as a
(two-sided) derivation with respect to the $\cdot $ product; then in
particular, $Sq_{1,1}(u;u)=Sq_{1}u$ for all $u\in H.$ Define $%
Sq_{p,q}=E_{p,q}:H^{\otimes p}\otimes H^{\otimes q}\rightarrow H$ for $%
p+q\geq 3$ by means of (\ref{dif}). Note that if the multiplicative
structure on $H$ is not free, such an extension might not exist. This
procedure gives a Hirsch algebra structure $\left\{ Sq_{p,q}\right\} $ on
the cohomology algebra $H$ in the following situations:

\begin{enumerate}
\item[\textit{(i)}] $H$ has trivial multiplication (e.g. the cohomology of a
suspension).

\item[\textit{(ii)}] $H$ is a polynomial algebra.

\item[\textit{(iii)}] $H$ has the following property: If $a\cdot b=0,$ then
\ $Sq_{1}a\cdot b=0=Sq_{1}a\cdot Sq_{1}b$ for all $a,b\in H.$
\end{enumerate}
\end{example}

Obviously we have the following proposition:

\begin{proposition}
\label{comparison} A morphism $f: A\rightarrow A^{\prime }$ of Hirsch
algebras induces a Hopf dga map of the bar constructions
\begin{equation*}
Bf : B A\rightarrow B A^{\prime }.
\end{equation*}
If the modules $A,A^{\prime}$ are $\Bbbk$-free and $f$ is a homology
isomorphism, so is $Bf.$
\end{proposition}

This proposition is useful when applying special models for a Hirsch algebra
$A$ to calculate the cohomology algebra $H^{\ast }(BA)=Tor^{A}({\Bbbk },{
\Bbbk })$ (see Subsection \ref{hohsec} below), and consequently, the loop
space cohomology $H^{\ast }(\Omega X;{\Bbbk })$ when $A=C^{\ast }(X;{\Bbbk }
) $ (see, for example, \cite{sanePOL}).

Given a Hirsch algebra $A$ with cohomology $H=H(A),$ let us construct a
Hirsch algebra model of $A.$ The commutative algebra $H$ admits a special \emph{
multiplicative} resolution $(RH,d),$ which is endowed with the Hirsch
algebra structure $\left\{ E_{p,q}\right\} .$ The perturbed differential $
d_{h}$ on $RH$ gives the desired Hirsch algebra model $(RH,d_{h})$ of $A.$

\subsection{Hirsch resolution}

Let ${H^{\ast }}$ be a graded algebra and recall that  a multiplicative
resolution $(R^{\ast }H^{\ast },d)$ of $H^{\ast }$ is the bigraded tensor
algebra $T(V)$ generated by the bigraded free $\Bbbk $-module
\begin{equation*}
V=\bigoplus_{j,m\geq 0}V^{-j,m},
\end{equation*}%
where $V^{-j,m}\subset R^{-j}H^{m}.$ The total degree of $R^{-j}H^{m}$ is the
sum $-j+m,$  $d$ is of bidegree $(1,0)$ and   $\rho :(RH,d)\rightarrow H$ is
a map of bigraded algebras
 inducing an isomorphism $\rho ^{\ast }:H^{\ast }(RH,d)
\overset{\approx }{\rightarrow }H^{\ast }$ where
$H^*$ is  bigraded via $H^{0,\ast }=H^{\ast}$
and $ H^{<0,\ast}=0$    (\cite{sane}; compare \cite{hal-sta}, \cite{hueb-kade}).
  In other words,
\begin{equation*}
\left( (R^{\ast }H^{m},d)\overset{\rho }{\rightarrow }H^{m}\right) =(\cdots
\overset{d}{\rightarrow }R^{-2}H^{m}\overset{d}{\rightarrow }R^{-1}H^{m}
\overset{d}{\rightarrow }R^{0}H^{m}\overset{\rho }{\rightarrow }H^{m})
\end{equation*}
is a usual free (additive) resolution of the $\Bbbk $-module $H^{m}$ for each
$m,$ and there is a multiplication on the family $\{R^{\ast }H^{m}\}_{m\in {
\mathbb{Z}}},$ which is compatible with both $d$ and the bidegree. When each
$H^{m}$ is $\Bbbk $-free, $\Omega BH$ (the cobar-bar construction of $H$) is
an example of $RH$ with $V=BH$. In general, the multiplicative structure of
 $H^{\ast }$ gives rise to (additively) \emph{non-minimal}
 submodules $
(R^{\ast }H^{m},d)$ even for $H^m$  to be  $\Bbbk$-free or $H^m=0.$
The reason for this is that  a (multiplicative) relation in $H$ involving elements of degree $<m$  can produce
an element $a\in R^{-1}H^{k}$ with $k<m,$ say $m=kn,$ some $n\geq2,$ and since the multiplication on $R^{\ast}H^{\ast}$  respects the bidegree,  the non-zero element $a^{n},$ the $n^{th}$ power of $a,$ ultimately belongs to    $ R^{-n}H^{m},$
the $n^{th}$ component of a $\Bbbk $-module resolution of $H^{m}$  (see the proof of Proposition \ref{hirschreso} below). Furthermore, even for $H$  to be a free cga over a field $\Bbbk,$  the
non-commutative nature of $RH$  fails to imply
 $R^*H^m$   to be a minimal $\Bbbk$-module resolution of $H^m,$     i.e.,
\begin{equation*}
R^{0}H^{m}=H^{m}\ \ \text{and}\ \ R^{-i}H^{m}=0,\,i>0;
\end{equation*}
 this is
quite different from the situation in \cite{hal-sta}.

For example, consider
the polynomial algebra $H={\mathbb{Z}}_{2}[x,y]$ with $x,y\in H^{2}$ and $
x_{0},y_{0}\in R^{0}H^{2}$ satisfying $\rho x_{0}=x$ and $\rho y_{0}=y.$
Then $R^{-1}H^{4}\neq 0$ since there is an element $a\in R^{-1}H^{4}$ such
that $da=x_{0}y_{0}+y_{0}x_{0}.$    In particular, if $H$ is the cohomology of
a dga $A$ with a \emph{non-commutative} $\smile _{1}$-product (and perhaps
higher order operations $E_{p,q};$ cf. Examples \ref{top} and \ref{hoh}),
then the construction of a \emph{Hirsch algebra} model of $A$ using $RH$
requires
 to add another element  $b$  in $R^{-1}H^{4}$ with $db=x_{0}y_{0}+y_{0}x_{0}.$
 Then denote $a=x_{0}\smile _{1}y_{0}$ and  $b=y_{0}\smile _{1}x_{0}$
 respectively (see
Theorem \ref{filtered}). Furthermore, if $H^{\ast }$ is $1$-reduced and we
wish to have a $1$-reduced multiplicative resolution $RH$, we must restrict
the resolution length of $R^{\ast }H^{m}$ so that $R^{-i}H^{m}=0$ for $i\geq
m-1$ (e.g. $H^{m}$ is $\Bbbk $-free for all $m$ or $H^{2}$ is ${\Bbbk }$
-free and $\Bbbk $ is a principal ideal domain). This motivates the
following definition:

\begin{definition}
\label{hirschreso} Let $H^{\ast }$ be a cga. An absolute Hirsch resolution of $H$ is a
multiplicative resolution
\begin{equation*}
\rho :R^{\ast }H^{\ast }\rightarrow H^{\ast },\ \ \ RH=T(V),\ \ \ V=\langle {
\mathcal{V}}\rangle ,
\end{equation*}
endowed with the Hirsch algebra structural operations
\begin{equation*}
E_{p,q}:RH^{\otimes p}\otimes RH^{\otimes q}\rightarrow V\subset RH
\end{equation*}
such that $V$ is decomposed as $V^{\ast ,\ast }=\mathcal{E}^{\ast ,\ast
}\oplus {U}^{\ast ,\ast }$ in which $\mathcal{E}^{0,\ast }=0,$ $U^{0,\ast
}=V^{0,\ast }$ and $\mathcal{E}^{\ast ,\ast }=\underset{{p,q\geq 1}}{
\bigoplus }\,\mathcal{E}_{p,q}^{<0,\ast }$ is distinguished by an
isomorphism of modules
\begin{equation*}
E_{p,q}:\bigoplus_{\substack{ i_{(p)}+j_{(q)}=s  \\ k_{(p)}+\ell _{(q)}=t}}
\left( \underset{1\leq r\leq p}{\otimes }R^{i_{r}}H^{k_{r}}\bigotimes
\underset{1\leq n\leq q}{\otimes }R^{j_{n}}H^{{\ell }_{n}}\right) \overset{
\approx }{\longrightarrow }\mathcal{E}_{p,q}^{s-p-q+1\,,\,t}\subset V^{\ast
,\ast }
\end{equation*}
where $x_{(r)}=x_{1}+\cdots +x_{r}.$
\end{definition}

 Given a  Hirsch algebra $(A,\{E_{p,q}\},d),$  a submodule $J\subset A$ is a \emph{Hirsch ideal}
of $A$ if it is an ideal with $E_{p,q}(a_{1},...,a_{p};a_{p+1},...,a_{p+q})
\in J$ whenever $a_{i}\in J$ for some $i.$

\begin{definition} Let $\rho_a:(R_aH,d)\rightarrow H$ be an absolute Hirsch resolution and $J\subset R_aH$ be a Hirsch ideal
such that $d:J\rightarrow J$ and the quotient map $g :R_aH\rightarrow R_aH/J$ is a homology isomorphism.
 A Hirsch resolution of $H$ is the Hirsch algebra $RH=R_aH/J$ with a map
$\rho:RH\rightarrow H$
such that $\rho_a=\rho\circ g.$
\end{definition}
Thus an absolute Hirsch resolution is a Hirsch resolution by taking $J=0.$

\begin{proposition}
\label{reso}
Every cga $H^{\ast }$ has an (absolute) Hirsch resolution $\rho :R^{\ast
}H^{\ast }\rightarrow H^{\ast }.$
\end{proposition}
\begin{proof}
We build a Hirsch resolution of $H^{\ast }$ by induction on the resolution
degree. Let $\mathcal{H}^{\ast }\subset H^{\ast }$ be a set of
multiplicative generators. Denote $\mathcal{V}^{0,\ast }=\mathcal{H}^{\ast };
$ let $V^{0,\ast }=\langle \mathcal{V}^{0,\ast }\rangle $ be the free $\Bbbk
$-module span of $\mathcal{V}^{0,\ast }$ and form the free (tensor) graded
algebra $R^{0}H^{\ast }=T(V^{0,\ast }).$ Obviously, there is a dga
epimorphism $\rho ^{0}:(R^{0}H^{\ast },0)\rightarrow H^{\ast }.$
Inductively, given $n\geq 0,$ assume we have constructed   a
$\Bbbk$-module $R^{(-n) }H^{\ast }=\oplus _{0\leq r\leq n}R^{-r}H^{\ast }$
with a map
$\rho ^{(n)}:(R^{(-n)}H^{\ast },d)\rightarrow H^{\ast }$  with
  $\rho^{r}(R^{-r}H^{\ast })=0$  for $1\leq r\leq n,$
where
  $d: R^{-r}H^{\ast }\rightarrow R^{-r+1}H^{\ast }$ is a differential of bidegree $(1,0)$ defined for
  $ 1\leq r\leq n$
  and
  acyclic in resolution degrees $-r$ for $1 \leq r<n;$
  $ R^{-r}H^{\ast }$ is  a  component of bidegree
$(-r,\ast)$ of $T(V^{(-r),\ast })$ for $V^{(-r),\ast }=V^{0,\ast }\oplus \dotsb \oplus V^{-r,\ast },$ so that
\[
R^{-r}H^{\ast}=V^{-r,\ast }\oplus \mathcal{D}^{-r,\ast} =\mathcal{E}^{-r,\ast }\oplus U^{-r,\ast }\oplus \mathcal{D}^{-r,\ast}
\]
where $\mathcal{E}^{-r,\ast }=
\underset{p,q\geq 1}{\bigoplus }\mathcal{E}_{p,q}^{-r,\ast }$ and
$\mathcal{E}_{p,q}^{-r,\ast }$ spans the set of (formal) expressions
\linebreak
 $
E_{p,q}(a_{1},...,a_{p};b_{1},...,b_{q}),\,a_{j}\in R^{-i_{k}}H^{\ast},\,b_{\ell}\in R^{-j_{\ell}}H^{\ast },
\,r=i_{(p)}+j_{(q)}+p+q-1,$
while
$\mathcal{D}^{-r,\ast}$ is the module of decomposables  of bidegree $(-r,\ast)$ in  $T(V^{(-r),\ast });$
 $d$  is given by  formula (\ref{dif}) on $\mathcal{E}^{-r,\ast },$  while acts as a  derivation on  $\mathcal{D}^{-r,\ast}.$

Let  $\mathcal{E}^{-n-1,\ast }=\underset{p,q\geq 1}{\bigoplus }\mathcal{E}
_{p,q}^{-n-1,\ast }$ where $\mathcal{E}_{p,q}^{-n-1,\ast }$  spans  the set of expressions
\linebreak
 $
E_{p,q}(a_{1},...,a_{p};b_{1},...,b_{q}),\,a_{k}\in R^{-i_{k}}H^{\ast },$
$b_{\ell}\in R^{-j_{\ell}}H^{\ast },\,n+1=i_{(p)}+j_{(q)}+p+q-1,$ and let
$\mathcal{D}^{-n-1,\ast}$ be the module of decomposables  of bidegree $(-n-1,\ast)$ in
$T\left(V^{(-n),\ast }\oplus \mathcal{E}^{-n-1,\ast }\right);$
Define $d$   by  formula (\ref{dif}) on $\mathcal{E}^{-n-1,\ast }$  and  as a  derivation on  $\mathcal{D}^{-n-1,\ast}$ so that
\[\mathcal{E}^{-n-1,\ast }\oplus \mathcal{D}^{-n-1,\ast }\overset{d}{\rightarrow} R^{-n}H^{\ast}\overset{d}{\rightarrow }R^{-n+\!1}H^{\ast}.   \]
  Define   a free $\Bbbk $-module $U^{-n-1,\ast}$ and $d$ on it   to achieve acyclicity in resolution
degree $-n,$ i.e, denoting
$
V^{-n-1,\ast }=\mathcal{E}^{-n-1,\ast }\oplus U^{-n-1,\ast },
$
we obtain a partial resolution for each $m\in \mathbb Z$
\begin{equation*}
V^{-n-\!1,m}\oplus \mathcal{D}^{-n-\!1,m}\overset{d}{\rightarrow }R^{-n}H^{m}\overset{d}{\rightarrow }R^{-n+\!1}H^{m}\overset{d}{\rightarrow }\cdots \overset{d}{
\rightarrow }R^{-\!1}H^{m}\overset{d}{\rightarrow }R^{0}H^{m}\overset{\rho }{
\rightarrow }H^{m}.
\end{equation*}
Define $R^{-n-1}H^{\ast}= V^{-n-\!1,\ast}\oplus \mathcal{D}^{-n-\!1,\ast} $ and
 $\rho^{n+1}:R^{-n-1}H^{\ast }\rightarrow H^{\ast }$ to be trivial.
 This completes the inductive step.

Finally, set $R^{\ast }H^{\ast }=\oplus _{n}R^{(-n)}H^{\ast }$ with $V^{\ast
,\ast }=\langle \mathcal{V}^{\ast ,\ast }\rangle ,$ $\mathcal{E}^{\ast ,\ast
}=\oplus _{n}\mathcal{E}^{-n,\ast },$\thinspace\ $U^{\ast ,\ast }=\oplus
_{n}U^{-n,\ast },$ $\rho |_{R^{0}H^{\ast }}=\rho ^{0}$ and $\rho
|_{R^{-n}H^{\ast }}=0$ for $n>0$ to obtain the desired resolution map $\rho
:RH\rightarrow H$.
\end{proof}

Note that in a Hirsch resolution   $(RH,  \{E_{p,q}\}, d),$  we may have relations among $E_{p,q}$'s
(e.g.  $E_{p,q}=0$ for some  $p,q\geq 1$;  cf. Subsection \ref{small}).
For example, the Hirsch structure of   $RH$  is
\emph{associative}  if the product $\mu _{_{E}}$ on the bar construction $B(RH)$ is associative  and is equivalent to the equalities among
$E_{p,q}$'s as follows.

Given a Hirsch algebra $A$ and an arbitrary triple
\begin{equation*}
(\mathbf{a};\mathbf{b};\mathbf{c})=(a_{1},...,a_{k}\,;b_{1},...,b_{\ell
}\,;c_{1},...,c_{r}),\ \ a_{i},b_{j},c_{s}\in A,
\end{equation*}
denote
\begin{multline*}
\mathcal{R}_{k,\ell ,r}((\mathbf{a};\mathbf{b});\mathbf{c})=\sum_{\substack{
k_{(p)}=k;\ell _{(p)}=\ell  \\ 1\leq p\leq k+\ell }}(-1)^{\varepsilon
}E_{p,r}(E_{k_{1},\ell _{1}}(a_{1},...,a_{k_{1}};b_{1},...,b_{\ell _{1}}), \\
\ \ \ \ \ \ \ \ \ \ \ \ \ \ \ \ \ \ \ \ \ \ \ \ \ \ \ \ \ \ \ \ \ \ \
...,E_{k_{p},\ell _{p}}(a_{_{k-k_{p}+1}},...,a_{k};b_{_{\ell -\ell
_{p}+1}},...,b_{p})\,;c_{1},...,c_{r})
\end{multline*}
and
\begin{multline*}
\mathcal{R}_{k,\ell ,r}(\mathbf{a};(\mathbf{b};\mathbf{c}))=\sum_{\substack{
\ell _{(q)}=\ell ;r_{(q)}=r \\ 1\leq q\leq \ell +r}}(-1)^{\delta
}E_{k,q}(a_{1},...,a_{k};E_{{\ell }_{1},r_{1}}(b_{1},...,b_{{\ell }
_{1}};c_{1},...,c_{r_{1}}), \\
...,E_{{\ell }_{q},r_{q}}(b_{_{\ell -{\ell }_{q}+1}},...,b_{\ell
};c_{_{r-r_{q}+1}},...,c_{q})),
\end{multline*}%
where we use the convention that $E_{0,1}(-;a)=E_{1,0}(a;-)=a,$ $
E_{0,m}(-;a_{1},...,a_{m})=E_{m,0}(a_{1},...,a_{m};-)=0,m\geq 2,$ and $
x_{(n)}=x_{1}+\cdots +x_{n},$ while the signs $\varepsilon $ and $\delta $
are induced by permutations of symbols $a_{i},b_{j},c_{s}$ (cf.\thinspace
\cite{Voronov2}). Then the associativity of $A$ is equivalent to the equalities
\begin{equation}\label{assotriple}
\mathcal{R}_{k,\ell ,r}((\mathbf{a};\mathbf{b});\mathbf{c})=\mathcal{R}
_{k,\ell ,r}(\mathbf{a};(\mathbf{b};\mathbf{c})),\,k,\ell ,r\geq 1.
\end{equation}

 Now consider    the expression
\begin{equation}\label{minus}
\mathcal{R}_{k,\ell ,r}(\mathbf{a};(\mathbf{b};\mathbf{c}))-\mathcal{R}
_{k,\ell ,r}((\mathbf{a};\mathbf{b});\mathbf{c}) \ \in\   \mathcal{E}^{1-k-\ell -r\, ,\, \ast }
\end{equation}
in an absolute Hirsch resolution $RH.$
 We have that this expression belongs to
$ \mathcal{E}^{-2,\ast }  $ and is a cocycle for $(\mathbf{a};\mathbf{b};\mathbf{c})=(a;b;c),$ $a,b,c\in
R^{0}H$ (\,see (\ref{3ass}) and Fig.\thinspace 1 below in which the
boundaries of both hexagons are labelled by the 6 components of
 $d\mathcal{R}
_{1,1,1}({a};({b};{c}))=d\mathcal{R}_{1,1,1}(({a};{b});{c})$\,).
So  there is an element, denoted by  $s(\mathcal{R}_{1,1,1}\left({a};({b};{c}))\right)\in V^{-3,*}$
such that $ ds(\mathcal{R}_{1,1,1}\left({a};({b};{c}))\right)=
\mathcal{R}_{1,1,1}({a};({b};{c}))-\mathcal{R}_{1,1,1}(({a};{b});{c}).
$
In general,
 define elements $s(\mathcal{R}_{k,\ell ,r}\left(\mathbf{a};(\mathbf{b};\mathbf{c}))\right)\in V$ such that
\begin{multline*}
d s(\mathcal{R}_{k,\ell ,r}\left(\mathbf{a};(\mathbf{b};\mathbf{c}))\right)+\\
s(\mathcal{R}_{k,\ell ,r}\left(d\mathbf{a};(\mathbf{b};\mathbf{c}))\right)+
(-1)^{\varepsilon_1}s(\mathcal{R}_{k,\ell ,r}\left(\mathbf{a};(d\mathbf{b};\mathbf{c}))\right)+
(-1)^{\varepsilon_2}s(\mathcal{R}_{k,\ell ,r}\left(\mathbf{a};(\mathbf{b};d\mathbf{c}))\right)\\
=  \mathcal{R}_{k,\ell ,r}(\mathbf{a};(\mathbf{b};\mathbf{c}))-\mathcal{R}
_{k,\ell ,r}((\mathbf{a};\mathbf{b});\mathbf{c}) \\
\varepsilon_1=|\mathbf{a}|+k,\,
\varepsilon_2=|\mathbf{a}|+|\mathbf{b}|+k+\ell.
   \end{multline*}
Consequently, $RH=R_aH/J_{ass}$ is an associative Hirsch resolution, where $J_{ass}\subset R_aH$ is a Hirsch ideal generated by
\[\{\mathcal{R}_{k,\ell ,r}(\mathbf{a};(\mathbf{b};\mathbf{c}))-\mathcal{R}
_{k,\ell ,r}((\mathbf{a};\mathbf{b});\mathbf{c}),\,
s(\mathcal{R}_{k,\ell ,r}\left(\mathbf{a};(\mathbf{b};\mathbf{c}))\right)      \}   .\]

 In particular, for $(\mathbf{a};\mathbf{b};
\mathbf{c})=(a;b;c)$ the  associativity of a Hirsch resolution implies the following
\begin{proposition}\label{assocup-one}
For  $a,b,c\in RH,$ there is the equality
\begin{multline} \label{3ass}
(a\smile _{1}b)\smile
_{1}c+E_{2,1}(a,b;c)+(-1)^{(|a|+1)(|b|+1)}E_{2,1}(b,a;c)  \\
=a\smile _{1}(b\smile
_{1}c)+E_{1,2}(a;b,c)+(-1)^{(|b|+1)(|c|+1)}E_{1,2}(a;c,b).
\end{multline}
\end{proposition}

A Hirsch resolution $(RH,d)$ is \emph{minimal} if
\begin{equation*}
d(u)\in \mathcal{E}+\mathcal{D}+\kappa _{u}\!\cdot \!V\ \ \text{for}\ \text{
all}\ u\in U,
\end{equation*}%
where ${\mathcal{D}^{\ast ,\ast }}\subset R^{\ast }H^{\ast }$ denotes the
submodule of decomposables $RH^{+}\!\cdot RH^{+}$ ($RH^{+}$ denotes $RH$
modulo the unital component) and $\kappa _{u}\in \Bbbk $ is non-invertible.
For example, when $\Bbbk =\mathbb{Z}$ we have $\kappa _{u}\in {\mathbb{Z}}%
\setminus \{-1,1\};$ when $\Bbbk $ is a field we have $\kappa _{u}=0$ for
all $u$. Note that a minimal Hirsch resolution is \emph{not} minimal in the
category of dgas since the resolution differential does not send
multiplicative generators into $\mathcal{D}$ even when $\Bbbk $ is a field.
Furthermore, the  notion of minimality of $RH$ does not depend upon whether some operation $E_{p,q}$  is zero
(cf. Subsection \ref{small}).
On the other hand, in order to define a $\smile_2$-operation in a simple way on $RH$  we  have to consider a non-minimal Hirsch
resolution  in the next subsection.

Such a  flexibility of choice of $RH$ is due to the trivial Hirsch structure of  $H,$ and, in practice,  the choice is suggested by a
Hirsch algebra $A$ that realizes $H$ as the cohomology algebra.

\subsection{QHHA structures on Hirsch algebras}
\label{cup2sec}

First, note that one  can introduce a $\smile_2$-product  on  a Hirsch resolution that satisfies  (\ref{realcup2}).
However, such a QHHA structure on $RH$ in not always satisfactory,  and we shall  consider  a $\cup_2$-operation simultaneously for the
reasons explained below.
For an even dimensional $a,$ or for any $a$
whenever $\nu=2,$ we have that  $a\smile_1 a$  is cocycle for $da=0;$ hence, there is an element $x\in RH$ with $dx=a\smile_1 a.$
But we can not identify $x$   with $a\smile _2a$ because $d(a\smile_2 a)=0$ according to (\ref{realcup2}).
On the other hand, it is helpful to denote   $x:=a\cup_2 a$ since
 certain formulas are  conveniently expressed in terms of the binary operation  $\cup_2$ (see, for example,  Proposition \ref{toperations} or Remark \ref{lie}).
 Furthermore, we can  identify $a\cup_2a$ with $\frac{1}{2} a\smile_2 a$ for $|a|$ even and 2  invertible in $\Bbbk.$

By construction of a Hirsch resolution in Proposition \ref{reso},
 the definition  of $\smile_2$ mimics  that of $\smile_1.$
We start with the consideration of
the expression \[(-1)^{a}a\smile _{1}b
+(-1)^{(|a|+1)|b|}b\smile _{1}a  \, \in \mathcal{E}^{-1,\ast }\ \  \text{for}  \  \   a,b\in {\mathcal V}^{0,*}. \]
 It is a   cocycle
in $(RH,d),$ and   hence, must be   killed by a multiplicative generator;    denote this generator  by  $  a\smile_2 b\in {
U}^{-2,\ast }.$ Inductively, assume that  the right-hand side of (\ref{realcup2})
has been  defined as an element of ${U}^{-n+1,\ast }.$    Then  it is bounded by a multiplicative generator   $a\smile_2 b\in
{U}^{-n,\ast }.$
Thus, $a\smile_2b\in U$ for all $a,b\in RH.$
In particular, if $dx=0,$ then $d(x\smile _2x)=0$  or  $d( \frac{\nu}{2}\,  x\smile _2x)=0$   for  $|x|$ to be odd  or for both $|x|$ and
$\nu$ to be even  respectively in which case a multiplicative generator $y\in U$ with $dy=x\smile _2x$ is denoted by $x\cup _3x.$

Now define a $\cup_2$-operation      by
\begin{equation} \label{cupab}
a\cup_2 b  =
\left\{
  \begin{array}{llll}
    a\smile_2 b , &   a\neq   b ,  &  a,b \text{ are in a basis of}\  RH   \\
    0,            &    a=b,    &      |a|\ \text{and}\ \nu \ \text{are odd } ,
  \end{array}
\right.
\end{equation}
while,  otherwise,  define $a\cup_2 a\in U $ by
\begin{equation} \label{cup22}
d(a\cup _{2}a)=
\left\{\!\!
  \begin{array}{lll}
\ \ \ \, a\smile_1 a+ a \smile_2 da + da \cup_3 da  , &  |a|\  \text{is even} \\
 \frac{\nu}{2}(a\smile_1 a+ a \smile_2 da) + da \cup_3 da    ,&  |a|\  \text{is odd},&  \nu \   \text{is even}.
  \end{array}
\right.
\end{equation}
Hence,  $a\cup_2 b\in U$ for  any $a,b\in RH,$ and let
\begin{equation*}
\mathcal{T}=\{ a\cup _2 b \in U\,|\,a,b\in RH \}.
\end{equation*}
Thus, we obtain the decomposition $U= \mathcal{T}\oplus {\mathcal{M}}, $
   some ${\mathcal{M}},$ and, hence, the decomposition
\begin{equation*}
V=\mathcal{E}\oplus U=\mathcal{E}\oplus \mathcal{T}\oplus {\mathcal{M}}.
\end{equation*}
In particular,  $\mathcal{T}$ contains elements of the form
 $a_{1}\cup _{2}\cdots \cup _{2}a_{n},$ $a_i\in RH,$
obtained by the iteration of the $\cup_2$-product  for  $n\geq 2.$ In particular, for
  $a_i\in V^{0, 2r}$     we have
 the following equality
 \begin{equation*}
d(a_{1}\cup _{2}\cdots \cup _{2}a_{n})=\sum_{(\mathbf{i};\mathbf{j})}sgn(
\mathbf{i};\mathbf{j})(a_{i_{1}}\cup _{2}\cdots \cup
_{2}\,a_{i_{k}})\,\smile _{1}\,(a_{j_{1}}\cup _{2}\cdots \cup
_{2}\,a_{j_{\ell }}),
\end{equation*}
where the summation is over unshuffles $(\mathbf{i};\mathbf{j}
)=(i_{1}<\cdots <i_{k}\,;j_{1}<\cdots <j_{\ell })$ of $\underline{n}$ with $
(a_{i_{1}},...,a_{i_{k}})=(a_{i_{1}^{\prime }},...,a_{i_{k}^{\prime }})$ if
and only if $\mathbf{i}=\mathbf{i}^{\prime }$ and $sgn(\mathbf{i};\mathbf{j})
$ is induced by the permutation sign $a_{i}\cup
_{2}a_{j}=(-1)^{|a_{i}||a_{j}|}a_{j}\cup _{2}a_{i}$ (see also Fig.\thinspace
1 for $n=3$); consequently, for $a_{1}=\cdots =a_{n}=a$ and $a^{\cup _{2}n} :=a\cup _{2}\cdots \cup _{2}a,$ we
get
\begin{equation}
da^{\cup _{2}n}=\sum_{k+\ell =n}a^{\cup _{2}k}\smile _{1}\,a^{\cup _{2}\ell
},\ \ \ \ \ \ \ \ k,\ell \geq 1.  \label{cup2power}
\end{equation}
Note that the above equalities do not depend on the parity of $a_{i}$'s
when $\nu =2.$

\begin{remark}\label{remark1} 1. The definition of ${\mathcal T}$ does not depend on the  (Hirsch)  associativity of $RH.$

2.  In a minimal Hirsch resolution one can also minimize the module ${\mathcal T}$ as
\begin{equation*}
\mathcal{T}=\{ a\cup _2 b \in U\,|\,a,b\in {\mathcal M} \},
\end{equation*}
while  $a\cup _2 b$  for  $a,b\in RH$  is extended by certain \emph{derivation} formulas.   These  formulas are rather complicated, but
they  could be written down if necessary.

3.  The module ${\mathcal M}$ reflects the complexity    of the
multiplicative relations of the commutative algebra $H.$
\end{remark}
For example, if $H$ is a polynomial algebra and $RH$ is a minimal Hirsch resolution,
then ${\mathcal{M}}=
{\mathcal{M}}^{0,\ast }=V^{0,\ast }$ and, consequently,
 $RH$ is completely determined by the $\smile _{1}$- and $\cup
_{2}$-operations \cite{sanePOL} (see also Theorem \ref{freehoh} below).

\subsection{Some canonical syzygies in the
Hirsch resolution}

Below we give  topological interpretation of some canonical syzygies in the
Hirsch resolution $RH.$ In particular these syzygies reflect the non-associativity of the $\smile_1$-product. Remark  that
higher order   canonical syzygies should be also related with the combinatorics of permutahedra. In practice, such  relations are helpful to construct  small Hirsch resolutions $RH$ (cf. \cite{sanePOL}, see also Remark \ref{remark1} above).

\newpage

\unitlength=1.00mm \special{em:linewidth 0.4pt} \linethickness{0.4pt}
\begin{picture}(111.66,48.77)
\ \ \ \ \ \ \ \ \ \ \ \ \ \ \ \ \ \ \ \ \ \ \put(12.33,45.33){\line(1,0){30.33}}
\put(42.66,45.33){\line(0,-1){30.00}} \put(42.66,15.33){\line(-1,0){30.33}}
\put(12.33,15.33){\line(2,1){30.33}} \put(12.33,30.33){\line(2,1){30.33}}
\put(42.66,45.33){\circle*{1.33}} \put(12.33,45.33){\circle*{1.33}}
\put(12.33,30.66){\circle*{1.33}} \put(12.33,15.33){\circle*{1.33}}
\put(42.66,15.33){\circle*{1.33}} \put(42.66,30.33){\circle*{1.33}}

\put(27.99,30.33){\makebox(0,0)[cc]{$a \smallsmile\!_1\, (b\!\smallsmile\!_1\,c)$}} \
\ \ \ \put(27.33,11.10){\makebox(0,0)[cc]{${(a\!\smallsmile\!_1b)c}$}}
\put(33.66,20.00){\makebox(0,0)[cc]{$E_{12}(a;b,c)$}}
\put(21.13,41.00){\makebox(0,0)[cc]{$E_{12}(a;c,\!b)$}} \ \

\put(5.00,38.00){\makebox(0,0)[cc]{$(a\!\smallsmile\!_1c)b$}} \
\put(5.00,23.00){\makebox(0,0)[cc]{$a(b\!\smallsmile\!_1c)$}} \ \ \

\put(28.00,48.44){\makebox(0,0)[cc]{${c(a\!\smallsmile\!_1b)}$}}

\put(49.99,23.00){\makebox(0,0)[cc]{$b(a\!\smallsmile\!_1c)$}}
\put(49.99,38.00){\makebox(0,0)[cc]{$(b\!\smallsmile\!_1c)a$}}

\put(34.99,38.33){\makebox(0,0)[cc]{$_{a\smallsmile\!_1(cb)}$}}
\put(17.90,21.33){\makebox(0,0)[cc]{$_{a\smallsmile\!_1(bc)}$}}
\put(104.33,45.33){\line(0,-1){30.00}} \put(104.33,15.33){\line(-1,0){30.33}}
\put(74.00,15.33){\line(2,1){30.33}} \put(74.00,30.33){\line(2,1){30.33}}
\put(104.33,45.33){\circle*{1.33}} \put(74.00,45.33){\circle*{1.33}}
\put(74.00,30.66){\circle*{1.33}} \put(74.00,15.33){\circle*{1.33}}
\put(104.33,15.33){\circle*{1.33}} \put(104.33,30.33){\circle*{1.33}}

\put(89.66,30.33){\makebox(0,0)[cc]{$(a\!\smallsmile\!_1\,b) \!\smallsmile\!_1\,c$}}

\put(95.33,20.00){\makebox(0,0)[cc]{$E_{21}(a,b;c)$}}
\put(83.00,41.00){\makebox(0,0)[cc]{$E_{21}(b,a;c)$}}

\put(89.00,11.10){\makebox(0,0)[cc]{${a(b\!\smallsmile\!_1c)}$}}

\put(67.00,23.00){\makebox(0,0)[cc]{$(a\!\smallsmile\!_1b)c$}}

 \put(111.66,23.00){\makebox(0,0)[cc]{$(a\!\smallsmile\!_1c)b$}}

 \put(67.00,38.00){\makebox(0,0)[cc]{$b(a\!\smallsmile\!_1c)$}}

\put(111.66,38.00){\makebox(0,0)[cc]{$c(a\!\smallsmile\!_1b)$}}

\put(96.66,38.33){\makebox(0,0)[cc]{$_{(ba)\smallsmile\!_1c}$}}

\put(79.95,21.33){\makebox(0,0)[cc]{$_{(ab)\smallsmile\!_1c}$}}

\put(59.00,29.67){\makebox(0,0)[cc]{$=$}} \put(12.33,15.33){\line(0,1){30.00}}

\put(74.00,15.33){\line(0,1){30.00}} \put(74.00,45.33){\line(1,0){30.33}}

\put(89.00,48.77){\makebox(0,0)[cc]{${(b\!\smallsmile\!_1c)a}$}}

\end{picture}

\unitlength=1.00mm \special{em:linewidth 0.4pt} \linethickness{0.4pt}
\begin{picture}(108.33,40.33)
\ \put(85.00,16.00) \ \ \bezier{52}(73.33,20.66)(68.66,15.33)(73.66,11.33)
\bezier{48}(73.66,11.33)(77.66,16.66)(73.66,20.66)
\put(95.00,16.00){\oval(10.00,9.33)[]}
\bezier{52}(95.00,20.66)(90.33,15.33)(95.33,11.33)
\bezier{48}(95.33,11.33)(99.33,16.66)(95.33,20.66)
\put(95.33,25.66){\oval(10.00,9.33)[]}
\bezier{52}(95.33,30.32)(90.66,24.99)(95.66,20.99)
\bezier{48}(95.66,20.99)(99.66,26.32)(95.66,30.32) \put(94.99,21.00){\circle*{1.33}}
\bezier{100}(73.66,11.00)(83.99,17.33)(94.99,11.00)
\bezier{100}(73.33,11.00)(85.66,17.33)(94.99,11.00)
\bezier{88}(73.33,11.00)(82.99,14.33)(94.66,11.00)
\bezier{88}(73.66,11.00)(85.99,8.66)(94.99,11.33) \put(95.66,30.33){\circle*{1.33}}
\put(73.66,11.33){\circle*{1.33}} \put(94.66,11.33){\circle*{1.33}}
\put(73.33,25.66){\oval(10.00,9.33)[]}
\bezier{52}(73.33,30.32)(68.66,24.99)(73.66,20.99)
\bezier{48}(73.66,20.99)(77.66,26.32)(73.66,30.32) \put(73.66,21.00){\circle*{1.33}}
\put(73.66,30.33){\circle*{1.33}} \bezier{120}(73.66,30.33)(84.66,40.33)(95.66,30.33)
\bezier{120}(73.66,30.33)(84.33,40.33)(95.66,30.33)
\bezier{100}(73.66,30.33)(82.66,24.33)(95.66,30.33)
\bezier{100}(73.66,30.33)(82.66,24.33)(95.66,30.33)
\bezier{92}(73.66,30.33)(83.33,27.33)(95.66,30.33)
\bezier{92}(73.66,30.33)(84.99,33.66)(95.66,30.33)
\bezier{112}(73.66,11.33)(83.99,2.33)(94.66,11.33)
\bezier{112}(73.66,11.33)(82.66,2.33)(94.66,11.33) \

\put(84.00,21.33){\makebox(0,0)[cc]{${a\!\cup_2\!b\!\cup_2\!c}$}}

\put(84.33,4.00){\makebox(0,0)[cc]{$_{(a\cup_2b)\smallsmile\!_1c}$}} \ \ \ \

\put(84.33,38.33){\makebox(0,0)[cc]{$_{c\smallsmile\!_1(a\cup_2b)}$}} \
\put(60.66,16.00){\makebox(0,0)[cc]{$_{a\smallsmile\!_1(b\cup_2c)}$}} \
\put(60.66,26.00){\makebox(0,0)[cc]{$_{(a\cup_2c)\smallsmile\!_1b}$}} \
\put(108.33,25.33){\makebox(0,0)[cc]{$_{(b\cup_2c)\smallsmile\!_1a}$}} \
\put(108.33,15.33){\makebox(0,0)[cc]{$_{b\smallsmile\!_1(a\cup_2c)}$}}
\put(5.67,20.33){\line(1,0){16.67}} \put(6.34,20.33){\circle*{1.33}}
\put(22.34,20.33){\circle*{1.33}} \put(40.00,20.16){\oval(11.33,11.67)[]}
\put(34.33,19.66){\circle*{0.00}} \ \put(34.33,20.66){\circle*{1.33}}
\put(45.67,20.66){\circle*{1.33}} \ \put(6.67,17.00){\makebox(0,0)[cc]{$ab$}}
\put(22.34,17.00){\makebox(0,0)[cc]{$ba$}} \

\put(14.00,23.33){\makebox(0,0)[cc]{$a\smile_1b$}} \
\put(40.00,28.33){\makebox(0,0)[cc]{$_{b\smile_1a}$}} \
\put(40.00,20.33){\makebox(0,0)[cc]{$a\cup_2\!b$}} \
\put(40.00,12.00){\makebox(0,0)[cc]{$_{a\smile_1b}$}}
\put(73.33,16.00){\oval(10.00,9.33)[]}
\end{picture}


\begin{center}
Figure 1. Topological interpretation of some canonical syzygies in the
Hirsch resolution $RH.$
\end{center}

\vspace{5mm}

The symbol\ "$=$"\ in the figure above assumes equality (\ref{3ass}); the picture for $a\cup _{2}b\cup _{2}c$ is in fact 4-dimensional and must
be understood as follows: Whence $a\cup _{2}b$ corresponds to the 2-ball,
the boundary of $a\cup _{2}b\cup _{2}c$ consists of the six 3-balls each of
which is subdivided into four 3-cells by fixing two equators (these cells just
correspond to the four summand components of the differential evaluated on
the compositions of the $\smile _{1}$- and $\cup _{2}$-products). Then given
a 3-ball, two cells from these four cells are glued to the ones of the
boundary of the (diagonally) opposite 3-ball, and the other cells are glued
to the ones of the boundaries of the neighboring 3-balls according to the
relation
\begin{equation*}
x\smile _{1}(y\smile _{1}z)+(x\smile _{1}y)\smile _{1}z=y\smile _{1}(x\smile
_{1}z)+(y\smile _{1}x)\smile _{1}z.
\end{equation*}

\subsection{Filtered Hirsch model}

\label{filteredalgebras}Recall that a dga $(A^{\ast },d)$ is \emph{multialgebra} if it is bigraded $A^{n}=\underset{n=i+j}{\oplus}
A^{i,j}$ $,i\leq 0,$ $
j\geq 0,$ and $d=d^{0}+d^{1}+\dotsb +d^{n}+\dotsb $ with $
d^{n}:A^{p,q}\rightarrow A^{p+n,q-n+1}$ \cite{hueb}. A dga $A$ is bigraded
via $A^{0,\ast }=A^{\ast }$ and $A^{i,\ast }=0$ for $i\neq 0;$ consequently,
$A$ is a multialgebra. A multialgebra $A$ is \emph{homological} if $d^{0}=0$
(hence $d^{1}d^{1}=0$) and
\begin{equation*}
H^{i}(\cdots \overset{d^{1}}{\rightarrow }A^{i,\ast }\overset{d^{1}}{
\rightarrow }A^{i+1,\ast }\overset{d^{1}}{\rightarrow }\cdots \overset{d^{1}}
{\rightarrow }A^{0,\ast })=0,\ \ i<0.
\end{equation*}%
For a homological multialgebra the sum $d^{2}+d^{3}+\dotsb +d^{n}+\dotsb $
is called a \emph{perturbation } of $d^{1}.$ In the sequel we always
consider homological multialgebras, $d^{1}$ is denoted by $d,$ $d^{r}$ is
denoted by $h^{r},$ and the sum $h^{2}+h^{3}+\dotsb +h^{n}+\dotsb $ is
denoted by $h.$ We sometimes denote $d+h$ by $d_{h}.$

A \emph{multialgebra morphism} $\zeta :A\rightarrow B$ between two
multialgebras $A$ and $B$ is a dga map of total degree zero that preserves
the    resolution (column)  filtration, so that $\zeta $ has the components $
\zeta =\zeta ^{0}+\dotsb +\zeta ^{i}+\dotsb ,\ \ \zeta
^{i}:A^{s,t}\rightarrow B^{s+i,t-i}.$ A chain homotopy $s:A\rightarrow B$
between two multiplicative maps $f,g:A\rightarrow B$ is an $(f,g)$-\emph{
derivation} homotopy if $s(ab)=s(a)g(b)+(-1)^{|a|}f(a)s(b).$ A \emph{homotopy
} between two morphisms $f,g:A\rightarrow B$ of multialgebras is an $(f,g)$
-derivation homotopy $s:A\rightarrow B$ of total degree $-1$ that lowers the
column filtration by 1.

A multialgebra is \emph{quasi-free} if it is a tensor algebra over a bigraded ${
\Bbbk }$-module. Given $m\geq 2,$ the map $h^{m}|_{A^{-m,\ast }}:A^{-m,\ast
}\rightarrow A^{0,\ast }$ is referred to as the \emph{transgressive}
component of $h$ and is denoted by $h^{tr}.$ A multialgebra $A$ with a
Hirsch algebra structure
\begin{equation*}
E_{p,q}:\otimes _{r=1}^{p}A^{i_{r},k_{r}}\bigotimes \otimes
_{n=1}^{q}A^{j_{k},{\ell }_{n}}{\longrightarrow }A^{s-p-q+1\,,\,t}
\end{equation*}%
with $(s,t)=\left(
i_{(p)}+j_{(q)}\,,\,k_{(p)}+\ell _{(q)}\right) ,\,p,q\geq 1,$ is called \emph{Hirsch
multialgebra}. A \emph{homotopy} between two morphisms $f,g:A\rightarrow
A^{\prime }$ of Hirsch (multi)algebras is a homotopy $s:A\rightarrow
A^{\prime }$ of underlying (multi)algebras and


\begin{equation}
\label{hirschhom}
\begin{array}{lll}
s(E_{p,q}(a_1,...,a_p\,;b_1,...,b_q))\\
\hspace{0.3in}=
 \underset{1\leq \ell\leq q}{\sum}
 (-1)^{\epsilon^a_{p}+\epsilon^b_{\ell-\!1}}
E_{p,q}(fa_1,\!...,fa_p\,;fb_1,\!...,fb_{\ell-1},sb_{\ell},gb_{\ell+1},\!...,
gb_q)$\vspace{1mm}$\\
\hspace{0.4in}
  + \!\!\underset{1\leq k\leq p}{\sum} (-1)^{\epsilon^a_{k-\!1}}
E_{p,q}(fa_1,...,fa_{k-1},sa_{k},ga_{k+1},...,ga_p\,;gb_1,...,gb_q)$\vspace{1mm}$\\
\hspace{0.4in}
 -\!\!\! \underset{\substack{1\leq i\leq p\\1<\ell\leq j\leq q}}{\sum}\!\!
 (-1)^{\epsilon_ {_{i\!,j\!,\ell}}}
E_{i,j}(fa_1,...,fa_{i}\,;fb_1,...,fb_{\ell-1},sb_{\ell},gb_{\ell+1},...,gb_j)\\

\hspace{1.6in} \cdot\,
E_{p-i,q-j}(fa_{i+1},...,fa_{p-1},sa_{p}\,;gb_{j+1},...,gb_{q})$\vspace{1mm}$\\
\hspace{0.4in}
 -\underset{\substack{0\leq i< k\leq  p\\1\leq j\leq q}}{\sum}\!\!
 (-1)^{\epsilon_{_{i\!,j\!,k}}}
E_{i,j}(  fa_1,...,fa_{i}  \,; sb_1,gb_2,...,gb_j)\\

\hspace{0.77in} \cdot\,
E_{p-i,q-j}(fa_{i+1},\!...,fa_{k-1},sa_{k},ga_{k+1},\!...,ga_{p}\,;gb_{j+1},\!...,gb_{q}),
$\vspace{1mm}$\\
\hspace{1.75in}\epsilon_{_{i,j,m}}=\epsilon^a_{p-\!1}+\epsilon^b_{m-\!1}+(\epsilon^a_{p}+
\epsilon^a_{i})\epsilon^b_{j},\, p,q\geq 1,
\end{array}
\end{equation}
in which the first equality is
\begin{equation*}
s(a\smile _{1}\,b)=(-1)^{|a|+1}fa\smile _{1}\,sb+sa\smile
_{1}\,gb-(-1)^{(|a|+1)(|b|+1)}sb\cdot sa.
\end{equation*}%
Denote the homotopy classes of morphisms between two Hirsch (multi)algebras
by $[-,-].$

\begin{definition}
\label{filtereddef} A quasi-free Hirsch homological multialgebra $
(A,\{E_{p,q}\},d+h)$ is a \textbf{filtered Hirsch algebra} if it has the
following additional properties:

\begin{enumerate}
\item[\textit{(i)}] In $A=T(V)$ a decomposition
\begin{equation*}
V^{\ast ,\ast }=\mathcal{E}^{\ast ,\ast }\oplus U^{\ast ,\ast }
\end{equation*}%
is fixed where $\mathcal{E}^{\ast ,\ast }=\underset{p,q\geq 1}{\bigoplus}\mathcal{E}_{p,q}^{<0\,,\ast}$ is distinguished by an
isomorphism of modules
\begin{equation*}
E_{p,q}:A^{\otimes p}\otimes A^{\otimes q}\overset{\approx }{\longrightarrow
}\mathcal{E}_{p,q}\subset V,\,\,\,p,q\geq 1;
\end{equation*}

\item[\textit{(ii)}] The restriction of the perturbation $h$ to $\mathcal{E}$
has no transgressive components $h^{tr},$ i.e., $h^{tr}|_{\mathcal{E}}=0.$
\end{enumerate}
\end{definition}

Given a Hirsch algebra $B,$ a \emph{filtered Hirsch model} for $B$ is a filtered
Hirsch algebra $A$ together with a Hirsch algebra map $A\rightarrow B$ that
induces an isomorphism on cohomology. Our next proposition, which is a
Adams-Hilton type of statement, exhibits a basic property of filtered Hirsch
algebras:

\begin{proposition}
\label{fibrant} Let $\zeta :B\rightarrow C$ be a map of (filtered)Hirsch
algebras that induces an isomorphism on cohomology. If $A$ is a  filtered
Hirsch algebra, there is a bijection of sets of homotopy classes of
(filtered) Hirsch algebra maps
\begin{equation*}
\zeta _{_{\#}}:[A,B]\overset{\approx }{\longrightarrow }[A,C].
\end{equation*}
\end{proposition}

\begin{proof}
Discarding Hirsch algebra structures, the proof goes by induction on the
resolution grading and is similar to that of Theorem 2.5 in \cite{hueb} (see
also \cite{saneDerived}). The Hirsch algebra structure serves to specify a
choice of homotopy $s$ on the multiplicative generators $\mathcal{E}\subset
V.$ When constructing a chain homotopy $s:A\rightarrow C$ between two
multiplicative maps $f,g:A\rightarrow C,$ we can choose an $s$ on $\mathcal{E
}^{i,\ast }$ that satisfies formula (\ref{hirschhom}) in each step of the
induction.
\end{proof}

The basic examples of a filtered Hirsch algebra are provided by the
following theorem, which states our main result on Hirsch algebras:

\begin{theorem}
\label{filtered} Let $H$ be a cga and let $\rho :(RH,d)\rightarrow H$ be an absolute
Hirsch resolution. Given a Hirsch algebra $A,$ assume there exists an
isomorphism $i_{A}:H\approx H(A,d).$ Then

\begin{enumerate}
\item[\textit{(i)}] Existence. There is a pair $(h,f)$ where $%
h:RH\rightarrow RH$ is a perturbation of the resolution differential $d$ on $
RH$ and
\begin{equation*}
f:(RH,d+h)\rightarrow A
\end{equation*}
is a filtered Hirsch model of $A$ such that $(f|_{_{R^{0}H}})^{\ast
}=i_{A}\rho |_{_{R^{0}H}}:R^{0}H\rightarrow H(A).$ \vspace{1mm}

\item[\textit{(ii)}] Uniqueness. If $(\bar{h},\bar{f})$ and  $\bar{f}:(RH,d+{
\bar{h}})\rightarrow A$ satisfy the conditions of (i), there is an
isomorphism of filtered Hirsch models
\begin{equation*}
\zeta :(RH,d+h)\overset{\approx }{\longrightarrow }(RH,d+{\bar{h}})
\end{equation*}%
of the form $\zeta =Id+\zeta ^{1}+\dotsb +\zeta ^{r}+\dotsb \ $with$\ \zeta
^{r}:R^{-s}H^{t}\rightarrow R^{-s+r}H^{t-r}$ such that $f$ is homotopic to $
\bar{f}\circ \zeta .$
\end{enumerate}
\end{theorem}
 Note that the proof of the theorem uses an induction on resolution grading  as it is used  by the  construction of  filtered model due to
  Halperin-Stasheff \cite{hal-sta} (compare also \cite{sane}, \cite{saneDerived}); although in the rational case for the existence and the uniqueness  of
  a pair $(h,f)$   the zero characteristic of $\Bbbk$ is essentially  involved, the proof below shows that such a restriction can be simply avoided. Here a technical subtlety  is that we have certain canonically chosen multiplicative generators on which $(h,f)$ must act by a canonical rule.
 \begin{proof}
\emph{Existence.} Let $RH=T(V)$ with $V=\mathcal{E}\oplus U.$ We define a
perturbation $h$ and a Hirsch algebra map $f:(RH,d+{h})\rightarrow (A,d)$ by
induction on resolution (column) grading.
 First consider $R^{0}H=T(V^{0,\ast})\, (=\! T(U^{0,\ast})).$ Define a chain map
${\mathfrak f}^{0}:(V^{0,\ast},0)\rightarrow (A,d)$   by  $({\mathfrak f}^{0})^{\ast }=i_{A}\rho
|_{_{V^{0,\ast}}}: V^{0,\ast} \rightarrow H(A).$
Extend ${\mathfrak f}^{0}$ multiplicatively to obtain a dga map  ${f}^{0}:R^0H\rightarrow A.$
There is a map ${\mathfrak f}^{1}:V^{-1,\ast}\rightarrow A^{\ast -1}$
with $f^{0}d|_{V^{-1,\ast}}=d{\mathfrak f}^{1};$ in particular, choose ${\mathfrak f}^{1}$
on $\mathcal{E}^{-1,\ast }\,(=\mathcal{E}_{1,1}^{-1,\ast })$  defined by
the formula ${\mathfrak f}^{1}(a\smile _{1}\,b)=f^{0}a\smile _{1}\,f^{0}b$    for $a,b\in
R^{0}H.$ Then extend  ${\mathfrak f}^{0}+{\mathfrak f}^{1}$ multiplicatively  to obtain a dga map
${\mathfrak f}^{(1)}_{\#}: T(V^{(-1),\ast})\rightarrow (A,d);$ then denote the restriction of ${\mathfrak f}^{(1)}_{\#}$ to $R^{(-1)}H$
by $f^{(1)}:(R^{(-1)}H,d)\rightarrow (A,d).$

Inductively, assume that a pair $(h^{(n)},f^{(n)})$ has been constructed
that satisfies the following conditions:

\begin{enumerate}
\item[(1)] $h^{(n)}=h^{2}+\dotsb +h^{n}$ is a derivation on $RH,$

\item[(2)] Equality (\ref{dif}) holds on $R^{(-n)}H$ for $d+h^{(n)}$ in
which
\begin{multline*}
h^{r}E_{p,q}(a_{1},...,a_{p}\,;b_{1},...,b_{q})=\sum_{i=1}^{p}(-1)^{\epsilon
_{i-\!1}^{a}}E_{p,q}(a_{1},...,h^{r}a_{i},...,a_{q}\,;b_{1},...,b_{q}) \\
+\sum_{j=1}^{q}(-1)^{\epsilon _{p}^{a}+\epsilon
_{j-\!1}^{b}}E_{p,q}(a_{1},...,a_{p}\,;b_{1},...,h^{r}b_{j},...,b_{q}),\
2\leq r\leq n,
\end{multline*}

\item[(3)] $dh^{n}+h^{n}d+\sum_{i+j=n+1}h^{i}h^{j}=0,$

\item[(4)] $f^{(n)}:R^{(-n)}H\rightarrow A$ is the restriction of a dga map  ${\mathfrak f}^{(n)}_{\#}:T\!\left(V^{(-n),\ast}\right)\rightarrow A $ to $R^{(-n)}H$ for $f^{(n)}= {f}^0+\dotsb +  {f}^n;$

\item[(5)] $f^{(n)}(d+h^{(n)})=df^{(n)}$ on $R^{(-n)}H,$ and

\item[(6)] $f^{(n)}$ is compatible with the maps $E_{p,q}$ on $\mathcal{E}%
^{(-n),\ast }.$
\end{enumerate}

Consider
\begin{equation*}
f^{(n)}(d+h^{(n)})|_{V^{-n-1,\ast }}:V^{-n-1,\ast }\rightarrow A^{\ast -n-1};
\end{equation*}%
clearly $df^{(n)}(d+h^{(n)})=0.$ Define a linear map $h^{n+1}:U^{-n-1,\ast
}\rightarrow R^{0}H^{\ast -n}$ with $\rho
h^{n+1}=i_{A}^{-1}[f^{(n)}(d+h^{(n)})]$ and extend $h^{n+1}$ on $RH$ as a
derivation (denoting by the same symbol) with
\begin{equation*}
dh^{n+1}+h^{n+1}d+\sum_{i+j=n+2}h^{i}h^{j}=0
\end{equation*}%
and
\begin{multline*}
h^{n+1}E_{p,q}(a_{1},...,a_{p}\,;b_{1},...,b_{q})=\sum_{i=1}^{p}(-1)^{%
\epsilon
_{i-\!1}^{a}}E_{p,q}(a_{1},...,h^{n+1}a_{i},...,a_{q}\,;b_{1},...,b_{q}) \\
+\sum_{j=1}^{q}(-1)^{\epsilon _{p}^{a}+\epsilon
_{j-\!1}^{b}}E_{p,q}(a_{1},...,a_{p}\,;b_{1},...,h^{n+1}b_{j},...,b_{q}).
\end{multline*}%
Then there is a map ${\mathfrak f}^{n+1}:V^{-n-1,\ast }\rightarrow A^{\ast -n-1}$ such
that  it is compatible with $E_{p,q}$ on $\mathcal{E}^{-n-1,\ast }$ and
\begin{equation*}
f^{(n)}(d+h^{(n+1)})|_{V^{-n-1,\ast }}=d{\mathfrak f}^{n+1}.
\end{equation*}
Extend  ${\mathfrak f^{(n+1)}}:={\mathfrak f}^0+\cdots +{\mathfrak f}^{n+1}$ multiplicatively
to obtain a dga map
 $  {\mathfrak f}^{(n+1)}_{\#}: T(V^{(-n-1),\ast})\rightarrow A ;$  the restriction of $ {\mathfrak f}^{(n+1)}_{\#}$
to $R^{(-n-1)}H$ denote by
\begin{equation*}
f^{(n+1)}:R^{(-n-1)}H\rightarrow A.
\end{equation*}
 Thus the construction of
the pair $\left( h^{(n+1)},f^{(n+1)}\right) $  completes the inductive
step. Finally, a perturbation $h=h^{2}+\dotsb +h^{n}+\dotsb $ and a Hirsch
algebra map $f$ such that $f=f^{0}+\dotsb +f^{n}+\dotsb $ are
obtained as desired.

\emph{Uniqueness.} Using Proposition \ref{fibrant} we construct a
multialgebra morphism
\begin{equation*}
\zeta :(RH,d+h)\rightarrow (RH,d+{\Bar{h}}),
\end{equation*}%
$\zeta =\zeta ^{0}+\zeta ^{1}+\dotsb ,$ with $\bar{f}\circ \zeta \simeq f;$
in addition, it is easy to choose $\zeta $ with $\zeta ^{0}=Id.$
\end{proof}

\subsection{Filtered model for a QHHA}

\label{model-qhha} Referring to Subsection \ref{cup2sec}, this section
considers the compatibility of the perturbation $h$ and the Hirsch map $f$
with the $\cup _{2}$-product of $RH$ in Theorem \ref{filtered}.
Even if $A$ is a QHHA in the theorem, it is  impossible to obtain a \emph{QHHA map} $f$ which commutes with $\cup
_{2}$-products because   the compatibility of parameters $q(-;-)$ under $f$ is
obstructed. When $A$ is a $\mathbb{Z}_{2}$-algebra, for example, the
obstruction is caused by the non-free action of $Sq_{1}$ on $H.$
However, when $q(-;-)=0$  for the $\cup_2$-operation in $A$ (cf. Example \ref{qhha}), one can refine the perturbation $h$ in Theorem
\ref{filtered}
as it is stated  in  Proposition \ref{toperations}  below (in particular,  item (i) of this proposition is  an essential detail
of the proof of the main result in \cite{saneGeodesic}).

Let $\texttt{T}\subset \mathcal{T}$ be a submodule defined by
\begin{equation*}
\texttt{T}=\langle   a\cup_2 b\in {\mathcal{T}}\mid a\neq b\
\text{in a basis  of}\ {\mathcal{M}} \rangle.
\end{equation*}

For $\nu =2,$ let $Sq_{1}:H^{m}(A)\rightarrow H^{2m-1}(A)$ be the map from
Example \ref{cartan}.

\begin{proposition}
\label{toperations}
Let $A$ be a QHHA with $\cup_2$-operation satisfying (\ref{realcup2})  (e.g. $A$ is a special Hirsch algebra from Example 2).  Then in the
filtered Hirsch model $f:(RH,d_{h})\rightarrow A$ given by Theorem \ref
{filtered}, the perturbation $h$ can be chosen such that

\begin{enumerate}
\item[\textit{(i)}] $h^{tr}|_{\texttt{T}}=0;$

\item[\textit{(ii)}]  Let $\nu =2.$  Then for $z_{i}=h^{tr}(a^{\cup _{2}2^{i}})$ with $a\in R^{0}H,$
\begin{equation*}
\rho z_{1}=Sq_{1}(\rho a)\ \ \text{and}\ \ h(a^{\cup _{2}2^{n}})=\underset{%
1\leq i<n}{\sum }z_{i}\!\cup _{2}a^{\cup _{2}(2^{n}-2i)}+z_{n}.
\end{equation*}
\end{enumerate}
\end{proposition}

\begin{proof}
(i) First, remark that any element  of $ \texttt{T}$ satisfies (\ref{realcup2}) (cf. (\ref{cupab})).
 Following the construction of a pair $(h,f)$ in the proof of Theorem
\ref{filtered}, define $f$ for $a\cup _{2}b\in \texttt{T}^{-2,*}$ with $
a,b\in {\mathcal{V}}^{0,\ast }$ by the formula
\begin{equation}\label{fcup2}
f(a\cup _{2}b)=
fa\cup _{2}fb.
\end{equation}
Since (\ref{realcup2}),  $f$ is chain with respect to the resolution differential $d$ of $RH,$  so we can
 set $h^{2}(a\cup _{2}b)=0.$ Inductively, assume that for $a\cup_2b\in \texttt{T}^{-r,*},$  $2\leq r<n,$
 the map $f$ is  defined by (\ref{fcup2}), while $h$ is defined by
 \begin{equation}\label{hcup2}
 h(a\cup_2b)=  ha\cup_2 b+(-1)^{|a|} a\cup_2 hb.   \end{equation}
Then for $a\cup_2b\in \texttt{T}^{-n,*}$ define $h$ again by (\ref{hcup2}).
Clearly,
 $fd_{h}(a\cup _{2}b)$ is a cocycle in $A$   and is bounded by $fa\cup _{2}fb.$
Therefore, we can define $f$ on $a\cup_2 b$ by (\ref{fcup2}). Consequently, we  set $h^{tr}(a\cup_2 b)=0$
as required.

(ii) Since $f$ is a Hirsch map, it commutes with $\smile _{1}$-products and
the first equality follows from the definition of $Sq_{1}.$ The verification
of the second equality follows immediately from (\ref{cup2power}).
\end{proof}

\begin{remark}
Whereas $Sq_{1}$ induces the product on $H(BA),$ the transgressive values $
z_{i}$ in item (ii) of Proposition \ref{toperations} are closely related with
the existence of the symmetric Massey products of the element $\sigma ^{\ast
}(\rho a)\in H(BA)$ for the suspension map $\sigma ^{\ast }:H^{\ast
}(A)\rightarrow H^{\ast -1}(BA)$ $($compare Theorem \ref{krainest} and
Remark \ref{lie} below$)$: When $\sigma ^{\ast }(\rho z_{k})=0$ for $k<i$ $($
e.g. $z_{k}\in {\mathcal{D}^{0,\ast }})$, the cohomology class $\sigma
^{\ast }(\rho z_{i})$ is automatically identified with the symmetric Massey
product $\langle \sigma ^{\ast }(\rho a)\rangle ^{2^{i}}.$
\end{remark}

Unlike Example \ref{top}, the Hirsch algebra $A$ provided by the following
example does not have a $\smile _{2}$-product. This fact allows us to lift
a combination $a\smile _{1}b\pm b\smile _{1}a$ for cocycles $a,b\in A$ to
the cohomology level as a non-trivial (binary) product (see also Subsection
\ref{hohsec}).

\begin{example}
\label{hoh} It is known that the Hochschild cochain complex $C^{\bullet
}(P;P)$ of an associative algebra $P$ admits an HGA structure $($\cite{kade2}
, \cite{GJ}$),$ which is a particular Hirsch algebra. Furthermore, whereas
the Hochschild cohomology $H=H\left( C^{\bullet }(P;P)\right) $ is a cga, $H$
is also endowed with the binary operation $x\ast y$ defined for $x=[a]$ and $
y=[b]$ by $x\ast y=[a\circ b-(-1)^{(|a|+1)(|b|+1)}b\circ a],$ where $\circ $
$(=\smile _{1})$ is  Gerstenhaber's operation on the Hochschild cochain
complex. The $\ast $ product on the Hochschild cohomology is referred to as
the G-algebra structure. Since $H$ is a cga, we can apply Theorem \ref%
{filtered} for $A=C^{\bullet }(P;P)$ and obtain the filtered Hirsch model $
f:(RH,d+h)\rightarrow C^{\bullet }(P;P).$ Given $a,b\in V^{0,\ast },$  obviously we
have $\rho h^{2}(a\cup _{2}\,b)=\rho a\ast \rho b$
(since $f^1(a\smile_1 b)=f^0a\circ f^0b$). In other words, the
non-triviality of the G-algebra structure on $H$ implies the non-triviality
of perturbation $h^{2}$ restricted to the submodule $\mathcal{T}\subset V.$
Consequently, the operation $a\cup _{2}\,b$ with $q(a,b)$ satisfying item $(
\ref{cup2})_{2}$ does not exist on the filtered Hirsch model of
 $C^{\bullet}(P;P)$ in general.
\end{example}

\subsection{A small Hirsch resolution $R_{\varsigma}H$}
\label{small}

 Let $A$ be a Hirsch algebra over $\Bbbk .$ Whereas
\linebreak
$
(RH,d_{h})=(T(V),d_{h})$ in a filtered Hirsch model $f:(RH,d_{h})\rightarrow
A,$ the calculation of $H(BA)$ can be carried out in terms of $V$ as
follows. Denote $\bar{V}=s^{-1}(V^{>0})\oplus \Bbbk $ and define the
differential $\bar{d}_{h}$ on $\bar{V}$ by the restriction of $d+h$ to $V$
to obtain the cochain complex $(\bar{V},\bar{d}_{h}).$ There are
isomorphisms
\begin{equation}
H^{\ast }(\bar{V},\bar{d}_{h})\approx H^{\ast }(B(RH),d_{_{B(RH)}})\overset{%
Bf^{\ast }}{\approx }H^{\ast }(BA,d_{_{BA}})\approx Tor^{A}(\Bbbk ;\Bbbk ).
\label{Vbar}
\end{equation}%
In particular, for $A=C^{\ast }(X;\Bbbk )$ with $X$ simply connected (cf.
Example \ref{top}),
\begin{equation*}
H^{\ast }(\bar{V},\bar{d}_{h})\approx H^{\ast }(BC^{\ast }(X;\Bbbk
),d_{_{BC}})\approx H^{\ast }(\Omega X;\Bbbk ).
\end{equation*}

\begin{remark}\label{remark3}
Note that the first isomorphism of (\ref{Vbar}) is a consequence of a
general fact about tensor algebras \cite{F-H-T}, while the second follows
from Proposition \ref{comparison}.
\end{remark}

Furthermore, to conveniently involve the multiplicative structure of (\ref
{Vbar}), one can reduce $V$ at the cost of $\mathcal{E}\subset V$ in the
manner we shall describe.
Let  $J_{\varsigma }\subset R_aH$ be the Hirsch ideal of an absolute Hirsch resolution $R_aH$  generated by
\begin{equation*}
\{E_{p,q}(a_{1},...,a_{p};a_{p+1},...,a_{p+q}),
\, d E_{p,q}(a_{1},...,a_{p};a_{p+1},...,a_{p+q})                  \,|\,p+q\geq 3\}
\end{equation*}
with
\[
                 \begin{array}{rlll}
a_{1},...,a_{p}\in R_aH,&   a_{p+1}\in {V}, & p\geq 1,\,q=1\\
                  a_{1},...,a_{p+q}\in R_aH ,&                  &p\geq1,\,  q>1.

                 \end{array}
               \]
Then
\begin{equation*}
R_{\varsigma }H=R_aH/J_{\varsigma }
\end{equation*}
is a Hirsch resolution of $H.$
Indeed, using (\ref{dif}) we see that $d:J_{\varsigma}\rightarrow J_{\varsigma }$  and  $H(J_{\varsigma},d)=0.$        Thus  $g_{\varsigma
}:(R_aH,d)\rightarrow (R_{\varsigma }H,d)$  is a homology isomorphism. We have an obvious
projection $\rho _{\varsigma }:(R_{\varsigma }H,d)\rightarrow H$ such that $
\rho =\rho _{\varsigma }\circ g_{\varsigma }.$ Consequently, $\rho_{\varsigma }$ is also a resolution map.
Furthermore,
we have  $h:J_{\varsigma
}\rightarrow J_{\varsigma }\ $so that $(R_{\varsigma }H,d_{h})$ is a Hirsch
algebra (in fact an HGA) and $g_{\varsigma }$ extends to a quasi-isomorphism
of  filtered  Hirsch  algebras
 \begin{equation}\label{gvarsigma}
  g_{\varsigma }:(R_aH,d_{h})\rightarrow (R_{\varsigma}H,d_{h}).
\end{equation}

Thus, the Hirsch (HGA) structure of $R_{\varsigma }H=T(V_{\varsigma })$ is
generated by the $\smile _{1}$-product and $(\ref{dif})$ is equivalent to
the following two equalities:

\begin{enumerate}
\item
\textit{The (left) Hirsch formula.} For $a,b,c\in R_{\varsigma }H:$
\begin{equation}
c\smile _{1}ab=(c\smile _{1}a)b+(-1)^{(|c|+1)|a|}a(c\smile _{1}b)
\label{hirsch1}
\end{equation}

\item
 \textit{The (right) generalized Hirsch formula.} For $a,b\in
R_{\varsigma }H$ and $c\in V_{\varsigma }$ with $d_{h}(c)=\sum c_{1}\cdots
c_{q},\,c_{i}\in V_{\varsigma }:$
\begin{equation} \label{hirsch2}
ab \smile _{1} \!c
 = \left\{ \!
\begin{array}{llll}
a(b\smile _{1}\!c)+(-1)^{|b|(|c|+1)}(a\smile _{1}\!c)\, b, &  q=1,\vspace{5mm}
&  &  \\
a(b\smile _{1}\!c)+(-1)^{|b|(|c|+1)}(a\smile _{1}\!c)\, b\vspace{1mm} &  &  &
\\
+\underset{1\leq i<j\leq q}{\sum }\! (-1)^{\varepsilon}\,c_{1}\cdots c_{i-1}(a\smile _{1}\!c_{i})\,c_{i+1}\vspace{1mm} &  &  &  \\
\cdots c_{j-1}(b\smile _{1}c_{j})\, c_{j+1}\cdots c_{q}, & q\geq 2,\vspace{1mm} &  &  \\
&  &  &
\end{array}
\right.
\end{equation}
where  $\varepsilon=(|a|+1)\left(\epsilon _{i-1}^{c}+i+1 \right)\, - \, (|b|+1)\left(\epsilon _{j-1}^{c}+j\right).$
\end{enumerate}
\begin{remark}
\label{remark2
}First, Formula $(\ref{hirsch2})$ can be thought of as a
generalization of Adams' formula for the $\smile _{1}$-product in the cobar
construction \cite[p.\,36]{adams-hopf1} from $q=2$ to any $q\geq 2.$ Second,
the usage of $R_{\varsigma }H$ shows that the multiplication $\mu _{E}^{\ast
}$ on $H^{\ast }(BA)\approx H^{\ast }(\bar{V}_{\varsigma },\bar{d}_{h})$ is
in fact determined only by the $\smile _{1}$-product on $V_{\varsigma }.$
\end{remark}

\vspace{0.2in}

\label{filteredpoly}
 \unitlength=1.00mm \special{em:linewidth 0.4pt}
\linethickness{0.4pt}
\begin{picture}(90.00,80.60)

 \put(80.00,69.90){\circle*{0.85}}

\put(80.00,59.67){\circle*{0.85}}
\put(80.00,43.25){\circle*{0.85}}
 \put(80.00,33.67){\circle*{0.85}}
\put(80.00,22.67){\circle*{0.85}}

\put(35.00,74.67){\makebox(0,0)[cc]{$+$}}

\put(26.50,75.00){\makebox(0,0)[cc]{$_{d}$}}

\put(29.60,68.63){\makebox(0,0)[cc]{$_{d}$}}

\put(45.57,68.63){\makebox(0,0)[cc]{$_{d}$}}

\put(44.00,75.00){\makebox(0,0)[cc]{$_{d}$}}

\put(55.00,61.10){\makebox(0,0)[cc]{$_{d}$}}

\put(64.83,41.07){\makebox(0,0)[cc]{$_{\mp 2x_0^2}$}}

\put(78.50,41.07){\makebox(0,0)[cc]{${_{x^2(=0)}}$}}

\put(90.00,70.00){\makebox(0,0)[ll]{$6m+3$}}
\put(90.00,59.67){\makebox(0,0)[ll]{$6m+2$}}
\put(90.00,43.34){\makebox(0,0)[lt]{$4m+2$}}

\put(65.00,57.00){\makebox(0,0)[cc]{${_{\mp 6h^2{x}_2}}$}}

\put(80.33,57.00){\makebox(0,0)[cc]{${{\langle x \rangle^3}}$}}

\put(58.70,80.60){\makebox(0,0)[cc] {$_{x_0^{\uplus 3}\,(=
x_0\smallsmile_1(x_0\smallsmile_1x_0)+2E_{12}(x_0;x_0,x_0))}$}}

\put(65.00,72.70){\makebox(0,0)[cc]
{${_{{\mp}3((x_0\smallsmile_1x_0)x_0+x_0(x_0\smallsmile_1x_0))}}$}}

 \put(19.50,72.70){\makebox(0,0)[cc]{${c}_2$}}

\put(34.33,71.80){\makebox(0,0)[cc]{$_{_{-6x_2+3u_2 }}$}}


 \put(49.67,56.67){\makebox(0,0)[cc]{$_{2h^2{\mathfrak b}_2}$}}

\put(79.50,14.33){\makebox(0,0)[cc]{$H^{\ast}$}}
\put(64.00,14.33){\makebox(0,0)[cc]{$R^{0}H^{\ast}$}}
\put(49.00,14.33){\makebox(0,0)[cc]{$R^{-1}H^{\ast}$}}
\put(34.00,14.33){\makebox(0,0)[cc]{$R^{-2}H^{\ast}$}}
\put(19.00,14.33){\makebox(0,0)[cc]{$R^{-3}H^{\ast}$}}

\put(16.40,69.90){\makebox(0,0)[cc]{$\cdots $}}
\put(16.40,59.40){\makebox(0,0)[cc]{$\cdots $}}
\put(16.40,33.33){\makebox(0,0)[cc]{$\cdots $}}
\put(16.40,42.93){\makebox(0,0)[cc]{$\cdots $}}
\put(16.40,22.53){\makebox(0,0)[cc]{$\cdots $}}

\put(35.40,43.33){\circle*{1.33}}
\put(20.80,43.33){\circle*{1.33}}

\put(35.40,22.53){\circle*{1.33}}
\put(20.80,22.53){\circle*{1.33}}
\put(50.00,22.53){\circle*{1.33}}

\put(49.67,48.33){\circle*{1.33}}

\put(64.20,69.90){\circle*{1.33}}

\put(64.50,59.67){\circle*{1.33}}
\put(50.00,59.67){\circle*{1.33}}
\put(35.00,59.40){\circle*{1.33}}
\put(20.50,59.40){\circle*{1.33}}

\put(50.00,33.40){\circle*{1.33}}
\put(35.40,33.40){\circle*{1.33}}
\put(20.90,33.40){\circle*{1.33}}

 \put(65.00,22.67){\circle*{1.33}}

\put(64.73,43.34){\circle*{1.33}}
 \put(49.67,43.34){\circle*{1.33}}

\put(65.00,33.67){\circle*{1.33}}

 \put(20.00,70.00){\circle*{1.33}}
\put(35.00,70.00){\circle*{1.33}}
\put(50.00,70.00){\circle*{1.33}}

\put(35.00,77.67){\circle*{1.33}}

\put(20.00,70.00){\vector(1,0){14.50}}

\put(49.00,70.00){\vector(1,0){14.50}}

\put(65.00,70.00){\vector(1,0){14.50}}

 \put(20.00,70.00){\vector(2,1){14.33}}

\put(35.00,77.63){\vector(2,-1){14.50}}

 \put(35.00,70.00){\vector(1,0){14.50}}

\put(19.80,70.00){\vector(3,-1){29.73}}

 \put(35.00,70.00){\vector(3,-1){29.20}}

\put(49.60,59.67){\vector(1,0){14.50}}

\put(35.00,59.67){\vector(1,0){14.50}}

\put(20.00,59.67){\vector(1,0){14.50}}

\put(64.67,59.67){\vector(1,0){14.67}}

\put(35.00,43.33){\vector(1,0){14.27}}

\put(20.40,43.33){\vector(1,0){14.27}}

 \put(50.00,43.34){\vector(1,0){13.83}}

\put(65.00,43.34){\vector(1,0){14.33}}

\put(65.00,33.67){\vector(1,0){14.33}}

\put(50.00,33.67){\vector(1,0){14.33}}
\put(35.00,33.67){\vector(1,0){14.33}}
\put(20.40,33.67){\vector(1,0){14.33}}

\put(68.50,14.03){\vector(1,0){7.00}}

 \put(65.00,22.67){\vector(1,0){14.33}}

 \put(50.00,22.67){\vector(1,0){14.33}}
 \put(35.00,22.67){\vector(1,0){14.33}}
  \put(20.40,22.67){\vector(1,0){14.33}}

\put(47.00,37.00){\makebox(0,0)[cc]{$h^2$}}

\put(31.50,46.00){\makebox(0,0)[cc]{${c_1}(=\!\mathfrak{b}_1\!)$}}

\put(64.67,30.00){\makebox(0,0)[cc]{$_{h^2{\mathfrak{b}}_1}$}}

\put(33.00,63.00){\makebox(0,0)[cc]{$h^2$}}
\put(55.17,66.00){\makebox(0,0)[cc]{$h^2$}}


\put(72.00,16.00){\makebox(0,0)[cc]{$\rho$}}

 \put(41.60,15.96){\makebox(0,0)[cc]{$_{d}$}}
\put(39.50,14.03){\vector(1,0){3.00}}
\put(26.60,15.96){\makebox(0,0)[cc]{$_{d}$}}
 \put(24.50,14.03){\vector(1,0){3.00}}
 \put(57.00,15.96){\makebox(0,0)[cc]{$_{d}$}}
  \put(55.00,14.03){\vector(1,0){3.00}}

\put(35.00,43.33){\vector(3,1){14.37}}
 \put(49.67,48.67){\vector(3,-1){14.33}}

 \put(35.00,43.33){\vector(3,-1){29.33}}


\put(58.67,51.00){\makebox(0,0)[cc]{$_{x_0^{\uplus 2}\,(=x_0\smallsmile_1x_0)}$}}

 \put(51.33,41.30){\makebox(0,0)[cc]{$_{2x_1}$}}

\put(49.67,46.00){\makebox(0,0)[cc]{$_{+}$}}

\put(72.33,24.00){\makebox(0,0)[cc]{$\rho$}}
\put(64.67,20.00){\makebox(0,0)[cc]{$x_0$}}

 \put(80.50,20.50){\makebox(0,0)[cc]{$x$}}

\put(90.00,33.67){\makebox(0,0)[lt]{$4m+1$}}
\put(90.00,22.67){\makebox(0,0)[lt]{$2m+1$}}
\put(40.67,47.00){\makebox(0,0)[cc]{$_{d}$}}

\put(60.67,47.00){\makebox(0,0)[cc]{$_{d}$}}

\put(44.33,44.63){\makebox(0,0)[cc]{$_{d}$}}
\put(55.67,44.63){\makebox(0,0)[cc]{$_{d}$}}
\end{picture}

\begin{center}
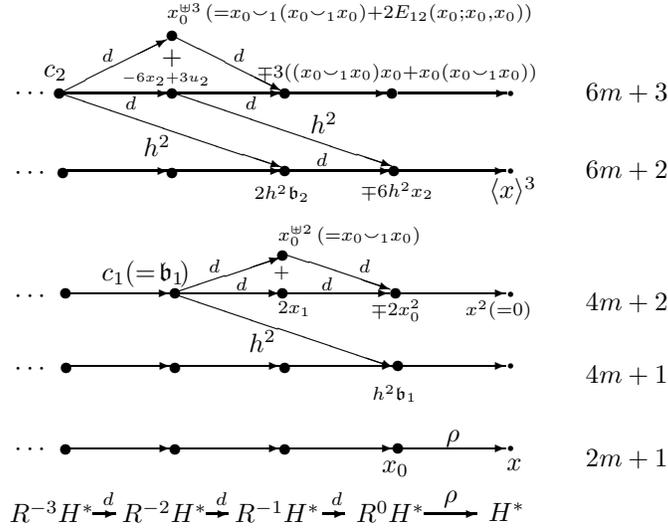
Figure 2. A fragment of the filtered Hirsch $\mathbb{Z}$-algebra obtained as
a perturbed resolution $(RH,d+h)$ of a cga $H.$\vspace{0.2in}
\end{center}

Note that for any Hirsch resolution of $H$ considered here, and consequently
for any filtered Hirsch model, the first two columns in Figure 2 are the
same.

\section{Some examples and applications}

In the discussion that follows we sometimes abuse notation and denote $
R_{\varsigma }H$ by $RH.$ As we mentioned in the introduction, certain
applications of the above material are given in \cite{sanePOL}, \cite
{saneBetti}. The applications that appear here are new.

\subsection{Symmetric Massey products}
\label{smss}
Recall  the definition of the $n$-fold symmetric
Massey product $\langle x\rangle ^{n}$ (cf. \cite{kraines}, \cite{may}).
Let $x\in H(A)$ be an element for a dga $A,$ and   $x_0\in A$ be a cocycle with $x=[x_0].$
Given $n\geq 3,$ consider a sequence $(x_0,x_1,...,x_{n-2})$ in $A$  such that
\begin{equation}\label{mass}
dx_{k}=\sum_{ i+j=k-1}
(-1)^{|x_{i}|+1}x_{i}x_{j},\ \ \ \    1\leq k\leq n-2;
\end{equation}
in particular, $dx_1=-(-1)^{|x_0|}x^2_0,$ i.e., $x^2=0.$ Then $\underset{i+j=n-2 }{\sum}
(-1)^{|x_{i}|+1}x_{i}x_{j}$
is a cocycle, and a subset of $H(A)$ formed by the classes of all such cocycles is denoted by
 $\langle x\rangle ^{n}.$
 (In other words, the existence of a sequence $(x_0,x_1,...,x_k,... )$  satisfying (\ref{mass}) for all $k$   implies that $c:=\underset{k\geq 0}\sum x_k   $ is a \emph{twisting} element in $A$  whenever this sum (possibly infinite) has a sense;
 an element $c\in A$  is twisting  if $dc=\pm\, c\cdot c;$ cf. \cite{berika3}.)

When $A=C^{\ast}(X; \mathbb{Z}_{p} )$ for  $p$ to be an odd prime, and
 $x\!\in \!H^{2m+1}(X;\mathbb{Z}_{p})$ is odd dimensional,
the following formula is established in \cite{kraines}
(for the dual case see \cite{kochman}):
\begin{equation}
{\langle x\rangle }^{p}=-\beta {\mathcal{P}}_{1}(x)  \label{kraines}
\end{equation}
where $\mathcal{P}_{1}:H^{2m+1}(X;\mathbb{Z}_{p})\rightarrow H^{2mp+1}(X;\mathbb{Z}_{p})$ is the Steenrod cohomology operation.
Thus,
 the formulas in \cite{kraines} and \cite{kochman} involve the
connection of the symmetric Massey products with the Steenrod and
Dyer-Lashof (co)homology operations in their respective topological
settings (cf. \cite{may}). Below Theorem \ref{krainest} emphasizes the algebraic
content of these formulas and generalizes them using a filtered Hirsch model over the integers.

\subsection{Massey syzygies  in the Hirsch resolution}

Let $(RH,d)$ be a Hirsch resolution of $H.$
Given a sequence of
relations of the form $da_{i}=\lambda b_{i}$ and
\begin{multline}
du_{i}=(-1)^{|a_{i}|+1}a_{i}a_{i+1}+\lambda v_{i},\ \ \
dv_{i}=(-1)^{|a_{i}|}b_{i}a_{i+1}+a_{i}b_{i+1},  \label{mprelation} \\
a_{i},u_{i},v_{i}\in RH,\ \lambda \in {\mathbb{Z}}\setminus \{-1,1\},\ \
1\leq i<n,
\end{multline}%
in $(RH,d),$ there are elements $u_{a_{i_{1}},...,a_{i_{k}}}\in RH,$ $3\leq
k\leq n,$ defined in terms of syzygies that mimic the definition of $k$-fold
Massey products arising from $k$-tuples $(a_{i_{1}},...,a_{i_{k}})$ \cite
{kraines}. Precisely, $u_{{a_{1},...,a_{n}}}$ is defined by
\begin{multline}
du_{{a_{1},...,a_{n}}}=\sum_{0\leq i<n}(-1)^{\epsilon _{i}^{a}}u_{{a_{1},...,a_{i}}}u_{{a_{i+1},...,a_{n}}}+
\lambda v_{{a_{1},...,a_{n}}},
\label{f-massey} \\
dv_{{a_{1},...,a_{n}}}=\sum_{0\leq i<n}((-1)^{\epsilon
_{i}^{a}+1}v_{a_{1},...,a_{i}}u_{a_{i+1},...,a_{n}}+u_{a_{1},...,a_{i}}v_{a_{i+1},...,a_{n}}),
\end{multline}%
with the convention that $u_{a_{i}}=a_{i},$ $u_{a_{i},a_{i+1}}=u_{i}$ and $
v_{a_{i}}=b_{i},$ $v_{a_{i},a_{i+1}}=v_{i}.$ When $b_{i}=0,$ equation (\ref
{f-massey}) reduces to
\begin{equation*}
du_{{a_{1},...,a_{n}}}=\sum_{0\leq i<n}(-1)^{\epsilon _{i}^{a}}u_{{
a_{1},...,a_{i}}}u_{{a_{i+1},...,a_{n}}}.
\end{equation*}%
We are interested in the special case of (\ref{mprelation}) obtained by
setting $a_{1}=\cdots =a_{n}.$ More precisely, we consider the following
situation (see also Example \ref{example} below).

Let $A$ be a torsion free Hirsch algebra over $\mathbb{Z}$ and fix a filtered model $
f:(RH,d_{h})\rightarrow A.$ For a module $C$ over $\mathbb{Z},$ let
$C_{\Bbbk }:=C\otimes _{\mathbb{Z}}\Bbbk $ and let $t_{\Bbbk }:C\rightarrow
C_{\Bbbk }$ be the standard map; then $A_{\Bbbk }=A\otimes _{\mathbb{Z}
}\Bbbk $ and $RH_{\Bbbk }=RH\otimes _{\mathbb{Z}}\Bbbk .$ Also let $H_{\Bbbk
}:=H(A_{\Bbbk }).$ There is the Hirsch model of $(A_{\Bbbk },d_{A_{\Bbbk }})$
given by
\begin{equation*}
f_{\Bbbk }=f\otimes 1:(RH_{\Bbbk },d_{h}\otimes 1)\rightarrow (A_{\Bbbk
},d_{A_{\Bbbk }}).
\end{equation*}
 Given an element $x\in H_{\Bbbk}$, let $x_{0}$
be a \emph{representative} of $x$ in $RH$ so that $[t_{_{\Bbbk}}f(x_{0})]=x.$ In particular, $x_{0}\in R^{0}H^{\ast }$ for $\beta
(x)=0,k\geq 1,$ and $x_{0}\in {R}^{-1}H^{\ast }$ with $
dx_{0}=\lambda x_{0}^{\prime },$ $x_{0}^{\prime }\in R^{0}H^{\ast },$ for $
\beta (x)\neq 0,$ where $\beta $ denotes the Bockstein cohomology
homomorphism associated with the sequence
\begin{equation*}
0\rightarrow {\mathbb{Z}}_{\lambda}\rightarrow {\mathbb{Z}}_{\lambda^{2}}\rightarrow {
\mathbb{Z}}_{\lambda}\rightarrow 0.
\end{equation*}
If $x\in H=H^{\ast }(A),$ then obviously $x_{0}\in R^{0}H^{\ast }.$ In any
case, assuming $x^{2}=0$ we have the corresponding relation in $(RH,d):$
\begin{equation*}
dx_{1}=(-1)^{|x_{0}|+1}x_{0}^{2}+\lambda x_{1}^{\prime }
\end{equation*}%
with the convention that $x_{1}^{\prime }=0$ whenever $x_{0}\in R^{0}H^{\ast
}.$ This equality is a special case of (\ref{mprelation}), so (\ref{f-massey}
) gives the following sequence of relations in $(RH,d):$
\begin{equation}
dx_{n}=\sum_{\substack{ i+j=n-1 \\ i,j\geq 0}}
(-1)^{|x_{i}|+1}x_{i}x_{j}+\lambda x_{n}^{\prime },\ \ \ \,n\geq 1,  \label{sms}
\end{equation}
where $x_{n}^{\prime }=0$ for $x_{0}\in R^{0}H.$

 We
have the following description of Massey symmetric products in terms of the sequence
${\mathbf{x}}=\{x_{n}\}_{n\geq 0}$ in $(RH,d_{h}).$ Denote $y_{i}=t_{\Bbbk
}x_{i}$ in $(RH_{\Bbbk },d_{h}).$ If $hy_{i}=0$ for $0\leq i<n,$ then (\ref
{sms}) implies $d_{h}d(y_{n})=dd(y_{n})=0,$ and consequently, $%
[dy_{n}]=-[hy_{n}].$ Therefore
\begin{equation} \label{massey}
f_{\Bbbk }^{\ast }[dy_{n}]=-f_{\Bbbk }^{\ast }[hy_{n}]\in \langle x\rangle
^{n+1}.
\end{equation}

Furthermore, the elements $
x_{n}$ appear in  a family of relations in $(RH,d).$ For example, these
relations can be deduced from the following observation. For $x\in H$ with $x^{2}=0,$
 let  $\iota :BH\rightarrow B(RH,d)$ be a chain map
such that
 $\iota ([\bar{x}|...|\bar{x}
])=(-1)^{n} [\, \overline {x_{n}} \,]$ for $[\bar{x}|...|\bar{x}]\in B^{n+1}H,$ $n\geq 0.$
 Assuming $BH$ is
endowed with the shuffle product $sh_{H},$ the map $\iota $ will be
multiplicative up to a chain homotopy $\mathfrak{b}.$ Since $B(RH)$ is
cofree, we can choose $\mathfrak{b}$ to be $(\, \mu _{E}\circ (\iota \otimes
\iota )\, , \iota \circ sh_{H})$-coderivation. This observation easily extends
to the $\!\!\!\mod \lambda $ case when $x_{0}\in R^{-1}H$ with $dx_{0}=\lambda x_{0}^{
\prime }.$ Now let \[\bar{\mathfrak{b}}_{k,\ell }:=
\mathfrak{b}(\,  \overset{k}{ \overbrace{[\bar{x}|...|\bar{x}]}}\otimes \overset{\ell }{\overbrace{[\bar{x}
|...|\bar{x}]}}\,)\,|_{_{\overline{RH}}}\ \
 \text{and}\ \
i_{[n]}:=i_{1}+\cdots +i_{n}+n;\]
then the equality $\mu _{E}(\iota \otimes \iota )-\iota \circ
sh_{H}=d_{_{B(RH)}}\mathfrak{b}+\mathfrak{b}d_{_{BH\otimes BH}}$ implies in
$(RH,d):$

For $|x_{0}|$ odd:
\begin{multline}
d\mathfrak{b}_{k,\ell }=(-1)^{k+\ell }\binom{k+\ell }{k}x_{k+\ell -1}
\label{odd} \\
\hspace{1.2in}+\sum_{\substack{ i_{[p]}=k,\,j_{[q]}=\ell }}(-1)^{k+\ell
+p+q}E_{p,q}(x_{i_{1}},...,x_{i_{p}};x_{j_{1}},...,x_{j_{q}}) \\
\hspace{0.7in}-\!\!\sum_{\substack{ _{0\leq r<k,0\leq m<\ell } \\
i_{[s]}=r,\,j_{[t]}=m}}(-1)^{r+m}\left(
(-1)^{s+t}\!E_{s,t}(x_{i_{1}},\!...,x_{i_{s}};x_{j_{1}},\!...,x_{j_{t}})
\mathfrak{b}_{k-r,\ell -m}\right.  \\
\left. +\,\binom{r+m}{r}\mathfrak{b}_{k-r,\ell -m}\,x_{r+m-1}\right) \!+
\lambda \mathfrak{b}_{k,\ell }^{\prime }
\end{multline}
in which $\mathfrak{b}_{k,\ell }^{\prime }=0$ for $x_{0}\in R^{0}H,$ and the
first equalities are:
\begin{equation*}
\begin{array}{lll}
d\mathfrak{b}_{1,1}=2x_{1}+x_{0}\smile _{1}x_{0}+\lambda \mathfrak{b}_{1,1}^{\prime
}, &  &  \\
d\mathfrak{b}_{2,1}=-3x_{2}+E_{2,1}(x_{0},x_{0};x_{0})-x_{1}\smile
_{1}x_{0}-x_{0}\mathfrak{b}_{1,1}+\mathfrak{b}_{1,1}x_{0}+\lambda \mathfrak{b}
_{2,1}^{\prime }, &  &  \\
d\mathfrak{b}_{1,2}=-3x_{2}+E_{1,2}(x_{0};x_{0},x_{0})-x_{0}\smile
_{1}x_{1}-x_{0}\mathfrak{b}_{1,1}+\mathfrak{b}_{1,1}x_{0}+\lambda \mathfrak{b}
_{1,2}^{\prime }. &  &
\end{array}
\end{equation*}
For $|x_{0}|$ even:
\begin{multline}
d\mathfrak{b}_{k,\ell }=(-1)^{k+\ell }\alpha _{k,\ell }\,x_{k+\ell -1}
\label{even} \\
\hspace{1.5in}+\sum_{\substack{ i_{[p]}=k,\,j_{[q]}=\ell }}(-1)^{k+\ell
+p+q}E_{p,q}(x_{i_{1}},...,x_{i_{p}};x_{j_{1}},...,x_{j_{q}}) \\
\hspace{0.7in}-\!\!\sum_{\substack{ _{0\leq r<k,0\leq m<\ell } \\
i_{[s]}=r,\,j_{[t]}=m}}\!\left(
(-1)^{(k+r+1)m+s+r+t}\!E_{s,t}(x_{i_{1}},\!...,x_{i_{s}};x_{j_{1}},
\!...,x_{j_{t}})\mathfrak{b}_{k-r,\ell -m}\right.  \\
\left. \hspace{2in}+\,(-1)^{k+\ell +r(\ell +m)}\alpha _{r,m}\mathfrak{b}
_{k-r,\ell -m}\,x_{r+m-1}\right) \!+\lambda \mathfrak{b}_{k,\ell }^{\prime }, \\
\\
\alpha _{i,j}=\left\{
\begin{array}{lll}
\binom{(i+j)/2}{i/2}, & i,j\ \text{are even},\vspace{1mm} &  \\
\binom{(i+j-1)/2}{i/2}, & i\ \text{is even},j\ \text{is odd},\vspace{1mm} &
\\
0, & i,j\ \text{are odd}, &  \\
&  &
\end{array}
\right.
\end{multline}
in which $\mathfrak{b}_{k,\ell }^{\prime }=0$ for $x_{0}\in R^{0}H,$ and the
first equalities are:
\begin{equation*}
\begin{array}{lll}
d\mathfrak{b}_{1,1}=x_{0}\smile _{1}x_{0}+\lambda \mathfrak{b}_{1,1}^{\prime }\ \ (
\text{i.e.,}\ \mathfrak{b}_{1,1}=x_{0}\cup _{2}x_{0}\ \ \text{when}\ \
x_{0}\in R^{0}H^{\ast }), &  &  \\
d\mathfrak{b}_{2,1}=-x_{2}+E_{2,1}(x_{0},x_{0};x_{0})-x_{1}\smile
_{1}x_{0}-x_{0}\mathfrak{b}_{1,1}-\mathfrak{b}_{1,1}x_{0}+\lambda \mathfrak{b}
_{2,1}^{\prime }, &  &  \\
d\mathfrak{b}_{1,2}=-x_{2}+E_{1,2}(x_{0};x_{0},x_{0})-x_{0}\smile
_{1}x_{1}+x_{0}\mathfrak{b}_{1,1}+\mathfrak{b}_{1,1}x_{0}+\lambda \mathfrak{b}
_{1,2}^{\prime }. &  &
\end{array}%
\end{equation*}%
Of course, for the sake of minimality, one can choose only certain $%
\mathfrak{b}_{k,\ell }$ above to be nontrivial. For example, let $|x|$ be
even, let ${\mathfrak{b}}_{2j+1}:={\mathfrak{b}}_{1,2j+1},$ and set $x_{2n}$
in (\ref{sms}) as
\begin{equation}
x_{2n}=-x_{0}\smile _{1}x_{2n-1}+\sum_{i+j=n-1}(x_{2i}{\mathfrak{b}}_{2j+1}-{
\mathfrak{b}}_{2j+1}x_{2i}).  \label{even1}
\end{equation}%
Thus one can also set ${\mathfrak{b}}_{1,2n}=0$ and eliminate $\mathfrak{b}
_{1,2n}$ from (\ref{even}); in particular, $\mathfrak{b}_{2,1}$ can be
identified with $x_{0}\smile _{2}x_{1}$ for $n=1.$

Note that for an HGA $A$ (e.g. $A=C^{\ast }(X;\mathbb{Z})$) we have that $
E_{p,q}=0$ for all $q\geq 2,$ that the second Hirsch formula up to homotopy
from Section \ref{Hiralgebras} becomes strict, and consequently, the
formulas above are much simpler (see also Subsection \ref{small}).
\begin{theorem}
\label{mzero}
Let $A$ be a Hirsch algebra over $\mathbb{Z}$ and let $\Bbbk $
be a field of characteristic $p\geq 0.$

\begin{enumerate}
\item[\textit{(i)}] Let $x\in H(A)$ with $x^{2}=0.$ If $\langle x\rangle ^{n}
$ is defined for $n\geq 3,$ it has a finite order.

\item[\textit{(ii)}] Let $x\in H_{\Bbbk }$ with $x^{2}=0$ and $p>0.$ Then $
\langle x\rangle ^{n}$ is defined for $3\leq n\leq p$ and vanishes whenever $
3\leq n<p.$

\item[\textit{(iii)}] Let $x\in H_{\Bbbk }$ with $x^{2}=0$ and $p=0.$ Then $
\langle x\rangle ^{n}$ is defined and vanishes for all $n.$
\end{enumerate}
\end{theorem}

\begin{proof}
(i) Observe that the inductive construction of the terms $h^{r},\,$\ $r\geq
2,$ of $h$ in $(RH,d_{h})$ implies $hx_{i}=0$ for $0\leq i\leq n-2$ whenever
$\langle x\rangle ^{n}$ is defined. Apply formulas (\ref{odd})--(\ref{even})
to verify that $m\langle x\rangle ^{n}=0$ with $m=n$ for $|x|$ odd (take $
(k,\ell )=(1,n-1)$ in (\ref{odd})), while $m=r-1$ or $m=r$ for $n=2r$ or $
n=2r+1$ (take $(k,\ell )=(2,n-2)$ in (\ref{even})) for $|x|$ even.

(ii)--(iii) The proof follows an argument similar to that in (i).
\end{proof}

\begin{remark}
First, regarding Theorem \ref{mzero}, item (i), note that formula $(\ref
{even1})$ implies that $\langle x\rangle ^{n}=0$ whenever $|x|$ and $n$ are
even$.$ Second, if $\left\vert x\right\vert $ is odd, formulas $(\ref{odd})$
--$(\ref{even})$ imply that whenever defined, $\langle x\rangle ^{n}$
consists of a single cohomology class independent of the parity of $n$ (see \cite{kraines}, \cite{kochman}).
\end{remark}

\subsection{The Kraines formula}
Let $p:=\lambda$ be an odd prime.
 Let $
a\in A^{2m+1}$ be an element with $da=0$ or $da=pa^{\prime }$ for some $
a^{\prime }.$ Given $n\geq 2,$ take (the right most) $n^{th}$-power of $\bar{
a}\in \bar{A}$ under the $\mu _{_{E}}$ product on $BA$ and consider its
component in $\bar{A}.$ Denote this component by $s^{-1}(a^{\uplus n})$ for $
a^{\uplus n}\in A^{2mn+1}.$ The element $a^{\uplus n}$ has the form
\begin{equation*}
a^{\uplus n}=a^{\smile _{1}n}+Q_{n}(a),
\end{equation*}
where $Q_{n}(a)$ is expressed in terms of $E_{1,k}$ for $1<k<n$ (for the
relations of small degrees involving this power, see also Fig. 2). For
example, $Q_{2}(a)=0\ $since $a^{\uplus 2}=a^{\smile _{1}2}$ and $
Q_{3}(a)=2E_{1,2}(a;a,a).$ In particular, if $A$ is an HGA, then obviously $
a^{\uplus n}=a^{\smile _{1}n}.$ Thus $da^{\uplus n}$ is divided by an
integer $p\geq 2$ if and only if $p$ is a prime and $n=p^{i},$ some $i\geq 1.
$ Consequently, the homomorphism
\begin{equation}\label{first}
\mathcal{P}_{1}:H_{\mathbb{Z}_{p}}^{2m+1}\rightarrow H_{\mathbb{Z}
_{p}}^{2mp+1},\ \ [t_{\mathbb{Z}_{p}}(a)]\rightarrow \lbrack t_{\mathbb{Z}
_{p}}(a^{\uplus p})],\ \ \ a\in A,\,\,d(t_{\mathbb{Z}_{p}}(a))=0,
\end{equation}
is well defined.

\begin{theorem}
\label{krainest} Let $A$ be a Hirsch algebra  as in Proposition \ref{toperations}. Let $A$ be  torsion free  and $p $ be an odd
prime. Then formula (\ref{kraines}) holds in $H_{\mathbb{Z}_p }$ for $
\mathcal{P}_1$ given by (\ref{first}).
\end{theorem}

\begin{proof}
Given $n\geq 1,$ let $\mathfrak{b}_{n}:=\mathfrak{b}_{1,n}$ and set $(k,\ell
)=(1,n)$ in (\ref{odd}) to obtain
\begin{multline}
d{\mathfrak{b}}_{n}=(-1)^{n+1}((n+1)x_{n}-\sum_{\substack{ j_{[q]}=n \\
1\leq q\leq n}}\!\!\!(-1)^{q}E_{1,q}(x_{0};x_{j_{1}},...,x_{j_{q}}))
\label{oddb} \\
\hspace{2in}+\sum_{i+j=n-1}\!\!(-1)^{i}\left( {\mathfrak{b}}_{j}x_{i}-x_{i}{
\mathfrak{b}}_{j}\right) +p{\mathfrak{b}}_{n}^{\prime }.
\end{multline}
By means of the element $x_{0}$ and the sequence $\{\mathfrak{b}
_{n}\}_{n\geq 1},$ form the sequence $\{c_{n}\}_{n\geq 1}$ in $RH$ as
follows:
\begin{equation*}
c_{1}=\mathfrak{b}_{1}\ \ \text{and}\ \ c_{n}=n!\,\mathfrak{b}
_{n}+x_{0}\smile _{1}c_{n-1},\,n\geq 2.
\end{equation*}
For $n=p-1,$ relation (\ref{oddb}) implies a relation of the form
\begin{equation}\label{cn}
dc_{p-1}=-p!\,x_{p-1}+x_{0}^{\uplus p}+pu_{p-1},
\end{equation}
where  $u_{p-1} \in   RH^+ \cdot RH^+ $  for $\beta(x)=0,$  while $u_{p-1} = w_{p-1} +(p-1)!\, {\mathfrak b}_{p-1}^{\prime}$ with
 $w_{p-1}\in   RH^+ \cdot RH^+ $  for $\beta(x)\neq 0.$   Hence, from $d^2(c_{p-1})=0$
 we get
\[d(x_{0}^{\uplus p})=p!\,dx_{p-1}-p\,du_{p-1} = p ((p-1)!\,dx_{p-1}-du_{p-1}).\]
 Obviously, $
h(x_{0}^{\uplus p})=0$ because       $h(x_0)=0$   (recall that a perturbation $h$  annihilates $R^{(-1)}H$ and  is
 a derivation on $\mathcal{E}$ ). Consequently,
 \[d_h(x_{0}^{\uplus p})= p((p-1)!\, dx_{p-1}-du_{p-1}).\]
 Taking into account  $(p-1)!=-1\!\!\mod p,$ and  passing to $H_{\mathbb{Z}_{p}}$ we obtain
\[
\beta {\mathcal{P}}_{1}(x)=
f_{_{{\mathbb{Z}}_{p}}}^{\ast}[- dy_{p-1}-dv_{p-1}   ]
=- f_{_{{\mathbb{Z}}_{p}}}^{\ast}[ dy_{p-1}]  - f_{_{{\mathbb{Z}}_{p}}}^{\ast} [dv_{p-1}]\ \text{for}\ v_{p-1}:=t_{_{{\mathbb{Z}}_{p}}}(u_{p-1}).
\]
Since $f_{_{{\mathbb{Z}}_{p}}}^{\ast}[ dy_{p-1}]={\langle x\rangle }^{p} $  by  (\ref{massey}), it remains to show that $f_{_{{\mathbb{Z}}_{p}}}^{\ast} [dv_{p-1}]=0.$
Indeed, if $\beta(x)=0,$ then $x_0\in R^0H,$ $u_{p-1}\in RH^+\cdot RH^+,$ and $hv_{p-1}=0$ by the similar argument as in the proof of Theorem  \ref{mzero} (ii). Consequently, $0=f_{_{{\mathbb{Z}}_{p}}}^{\ast} [-hv_{p-1}]= f_{_{{\mathbb{Z}}_{p}}}^{\ast} [dv_{p-1}].$
If $\beta(x)\neq 0,$ then $x_0\in R^{-1}H,$ and let $dx_0=px'_0.$ We have that $u_{p-1}$ contains $ \mathfrak{b}^{\prime}_{p-1}$ as a summand, and
$hv_{p-1}= -h \mathfrak{b}^{\prime}  _{p-1}.$
Denoting $z_0=g_{\varsigma}(x_0)$ and  $z'_0=g_{\varsigma}(x'_0)$
in  $(R_{\varsigma},d_h)$ where $g_{\varsigma}$ is given by (\ref{gvarsigma}),
 we have that $g_{\varsigma}(x_{0}^{\uplus p})=z_{0}^{\smile_1 p}$ and
$  g_{\varsigma}(h  \mathfrak{b}^{\prime} _{p-1})$ is mod p cohomologous to
\[\sum_{0\leq i<p} z_0^{\smile _{1} i}\smile _{1}z_0^{\prime}\smile _{1} z_0^{\smile_{1} p-i-1},\ \text{a summand component of}\ d(z_{0}^{\smile_1 p}).\]
 But this component   bounds
$\underset{0\leq i\leq p-2}{\sum} z_0^{\smile _{1} i}\smile _{1}(z_0\cup_2\, z_0^{\prime})\smile _{1} z_0^{\smile_{1} p-i-2}$ mod p
that finishes the proof.
\end{proof}

\begin{remark}
 When $p=2$ the
relation $d(x_{0}\smile _{1}x_{0})=-2x_{0}^{2}+2(x_{0}^{\prime }\smile
_{1}x_{0}+x_{0}\smile _{1}x_{0}^{\prime })$ implies the Adem relation $
Sq_{0}(a)=Sq^{1}Sq_{1}(a)$ in $H_{\mathbb{Z}_{2}}$ thought of as the
\textquotedblleft Kraines formula\textquotedblright\ ${\langle a\rangle }
^{2}=a^{2}=\beta Sq_{1}(a).$
\end{remark}

\begin{example}
\label{example} Fix a Hirsch filtered model $f:(RH,d_{h})\rightarrow A$ with
$RH=T(V).$ Suppose that we are given a single relation
\begin{equation}
da=\lambda b,\ \ a\in V^{-1,2k+1},b\in V^{0,2k+1},\,\lambda \geq 2,\,k\geq 1,
\label{single}
\end{equation}
and deduce the following relations in $(RH,d):$ First, define $c\in V$ by
\begin{equation}
dc=\left\{
\begin{array}{lll}
ab+\frac{\lambda }{2}b\smile _{1}b, & \lambda  & \text{is even} \\
2ab+\lambda b\smile _{1}b, & \lambda  & \text{is odd}.
\end{array}
\right.   \label{ab}
\end{equation}
When $\lambda $ is odd, denote $($cf. $(\ref{mprelation})$$)$
\begin{equation*}
u_{2a,b}:=-c,\ \ u_{b,2a}:=c-2a\smile _{1}b\ \ \text{and}\ \
u_{2b,b}:=2ab+\!(\lambda -1)b\smile _{1}b
\end{equation*}
 and obtain
\begin{equation}
\begin{array}{lll}
\ \ \,du_{a,a}=-a^{2}+\lambda v_{a,a},\ \ v_{a,a}=c-a\smile _{1}b, &  &  \\
du_{a,2b,b}=-au_{2b,b}-u_{a,2b}b+\lambda v_{a,2b,b} &  &  \\
\hspace{1.6in}=-2a^{2}b-(\lambda -1)a(b\smile _{1}b)+cb+\lambda u_{b,2b,b},
\vspace{1mm} &  &  \\
du_{b,2a,b}=bu_{2a,b}-u_{b,2a}b+\lambda v_{b,2a,b}=bc-(c-2a\smile
_{1}b)b+\lambda u_{b,2b,b}, &  &  \\
du_{a,2a,b}=-au_{2a,b}+u_{2a,a}b+\lambda v_{a,2a,b}, &  &
\end{array}
\label{system}
\end{equation}
where $v_{a,2b,b}=v_{b,2a,b}=u_{b,2b,b}=2u_{b,b,b}.$ Keeping in mind the
fact that $d_{h}^{2}=0,$ there is the following action of the perturbation $h
$ on the relations above:
\begin{equation*}
\begin{array}{lll}
\ \ \, dh^{2}u_{a,a}=-\lambda h^{2}c, &  &  \\
dh^{2}u_{a,2b,b}=-h^{2}c\cdot b-\lambda h^{2}u_{b,2b,b}, &  &  \\
dh^{2}u_{b,2a,b}=b\cdot h^{2}c+h^{2}c\cdot b-\lambda h^{2}u_{b,2b,b}, &  &
\\
dh^{2}u_{a,2a,b}=-a\cdot h^{2}c-2h^{2}u_{a,a}\cdot b-\lambda h^{2}v_{a,2a,b},
&  &  \\
dh^{3}u_{a,2a,b}=-h^{3}u_{2a,a}\cdot b-\lambda
h^{3}v_{a,2a,b}-h^{2}h^{2}u_{a,2a,b}. &  &
\end{array}%
\end{equation*}%
Below we shall exploit the third equality in list of relations above. First,
we have
\begin{equation*}
d\left( h^{2}u_{b,2a,b}+b\smile _{1}h^{2}c\right) =-\lambda h^{2}u_{b,2b,b}.
\end{equation*}%
Suppose that $\Bbbk $ is a ring such that $\nu $ divides $\lambda $ and
\begin{equation}
\begin{array}{lll}
\lbrack t_{\Bbbk }(a)][t_{\Bbbk }(b)]=0. &  &
\end{array}
\label{zero}
\end{equation}
By $(\ref{ab})$ one has $[t_{\Bbbk }(ab)]=-[t_{\Bbbk }h^{2}c],$ so that $
h^{2}c=0\!\!\!\ \mod \nu$ above. Denoting $[t_{\Bbbk }f(a)]:=y$ and $
[t_{\Bbbk }f(b)]:=x,$ we have $xy=0$ by $(\ref{zero}).$ Thus the triple
Massey product $\langle x,y,x\rangle $ is defined in $H_{\Bbbk }$ and
contains $[t_{\Bbbk }f(bu_{a,b}-u_{b,a}b)]$ $(= -[ t_{\Bbbk }f(hu_{b,a,b})] ).$
Obviously, $\langle x\rangle ^{3}$ is also defined and
\begin{equation*}
\beta _{\lambda }\langle x,y,x\rangle =-\langle x\rangle ^{3}
\end{equation*}
(here $\beta _{\lambda }$ denotes the Bockstein map associated with $
0\rightarrow {\mathbb{Z}}_{\nu }\rightarrow {\mathbb{Z}}_{\nu \lambda
}\rightarrow {\mathbb{Z}}_{\lambda }\rightarrow 0$). Now let $p=\lambda =3$
and consider $(\ref{cn})$ for $x.$ Then
\begin{equation*}
c_{2}=2\mathfrak{b}_{2}+x_{0}\smile _{1}\mathfrak{b}_{1},\ \ x_{0}^{\uplus
3}=x_{0}^{\smallsmile _{1}3}+2E_{1,2}(x_{0};x_{0},x_{0}),\ \ u_{2}=\mathfrak{
b}_{1}x_{0}-x_{0}\mathfrak{b}_{1}
\end{equation*}
and
\begin{equation*}
dc_{2}=-6x_{2}+x_{0}^{\smallsmile _{1}3}+2E_{1,2}(x_{0};x_{0},x_{0})+3(
\mathfrak{b}_{1}x_{0}-x_{0}\mathfrak{b}_{1}).
\end{equation*}%
Since $[x_{0}]^{2}=0,$ one has $h^{2}\mathfrak{b}_{1}=0$ and hence
\begin{equation*}
hc_{2}=2(h^{2}+h^{3}){\mathfrak{b}}_{2}
\end{equation*}
$($for the relations above, see also Fig. $2).$ In particular, $
dh^{2}c_{2}=6h^{2}x_{2}.$ Let $a:=y_{0},$ $b:=x_{0},$ $u_{b,b}:=x_{1}$ and $%
u_{b,b,b}:=x_{2}$ and set $h^{2}c_{2}=-2h^{2}u_{x_{0},y_{0},x_{0}}.$
Furthermore, if we also have $h^{3}c_{2}=h^{3}u_{x_{0},y_{0},x_{0}}\!\!\mod 3%
,$ then $[t_{\Bbbk }f(x_{0}^{\uplus 3})]=-[t_{\Bbbk }f(hc_{2})]=-[t_{\Bbbk
}f(hu_{x_{0},y_{0},x_{0}})]$ and, consequently,
\begin{equation}
{\mathcal{P}}_{1}(x)\in \langle x,y,x\rangle .  \label{powerm}
\end{equation}%
For example, let $A=C^{\ast }(BF_{4};{\mathbb{Z}}_{3}),$ the cochain complex
of the classifying space $BF_{4}$ of the exceptional group $F_{4}.$ Then
equality $(\ref{zero})$ together with $(\ref{powerm})$ holds in $H(BF_{4};{%
\mathbb{Z}}_{3}).$ More precisely, let $x_{i}\in H^{i}(BF_{4};{\mathbb{Z}}%
_{3})$ be multiplicative generators in notation of \cite{Toda} and recall
the following relations among them: $x_{8}x_{9}=0=x_{4}x_{21},$ $\delta
x_{8}=x_{9},$ $\delta x_{25}=x_{26};$ also $\mathcal{P}^{3}(x_{9})=x_{21}$
and $\mathcal{P}^{1}(x_{21})=x_{25};$ thus $\mathcal{P}^{1}\mathcal{P}%
^{3}(x_{9})=\mathcal{P}_{1}(x_{9})=x_{25}$ by an application of the Adem
relation. Thus the knowledge of both $H^{\ast }(BF_{4};{\mathbb{Z}}_{3})$
and $H^{\ast }(F_{4};{\mathbb{Z}}_{3})$ in low degrees enables us to use the
filtered Hirsch model of $BF_{4}$ to deduce the following: Let $a$ and $b$
be defined in $(\ref{single})$ by $[t_{_{\mathbb{Z}_{3}}}f(a)]=x_{8}$ and $%
[t_{_{\mathbb{Z}_{3}}}f(b)]=x_{9}.$ Then $[t_{_{\mathbb{Z}%
_{3}}}f(hc_{2})]=[t_{_{\mathbb{Z}_{3}}}f(hu_{b,a,b})]=-x_{25}$ and $[t_{_{%
\mathbb{Z}_{3}}}f(h^{2}u_{b,b,b})]=x_{26}$ so that
\begin{equation*}
\langle x_{9}\rangle ^{3}=-\beta {\mathcal{P}}_{1}(x_{9})\ \ \text{with}\ \ {%
\mathcal{P}}_{1}(x_{9})=\langle x_{9},x_{8},x_{9}\rangle .
\end{equation*}

Finally, we remark that the both sides of this formula become trivial under
the loop suspension map $\sigma ^{\ast }:H^{\ast }(BF_{4};{\mathbb{Z}}%
_{3})\rightarrow H^{\ast -1}(F_{4};{\ \mathbb{Z}}_{3})$ by a general
well-known fact about Massey products \cite{kraines}, \cite{kraines2}
(compare $\mathcal{P}_{1}(i_{3})$ for $i_{3}\in H^{3}(K({\mathbb{Z}}_{3};3);{%
\mathbb{Z}}_{3})).$
\end{example}

\subsection{Hochschild cohomology with the $G$-algebra structure}

\label{hohsec}In this section we assume that $\Bbbk $ is a field of
characteristic zero. Refer to Example \ref{hoh} and recall that the HGA
structure $E=\{E_{p,q}\}_{p\geq 0;q=0,1}$ on the Hochschild cochain complex $
\,A=C^{\bullet }(P;P)$ induces an associative product $\mu _{E}$ on the bar
construction $BA$ and hence the product $\mu _{E}^{\ast }$ on $H^{\ast
}(BA)=Tor_{\ast }^{A}(\Bbbk ,\Bbbk ).$ Since $Tor_{\ast }^{A}(\Bbbk ,\Bbbk )$
is an associative algebra, it can be converted into a Lie algebra in the
standard way.

\begin{theorem}
\label{freehoh} If the Hochschild cohomology $H^*=H(C^{\bullet}(P;P))$ is a
free algebra, then the Lie algebra structure on $Tor^A_*(\Bbbk,\Bbbk)$ is
completely determined by that of the $G$-algebra $H^*.$ Consequently, the
product $\mu^*_{E}$ on $Tor^A_*(\Bbbk,\Bbbk)$ is commutative if and only if
the $G$-product on $H^*$ is trivial.
\end{theorem}

\begin{proof}
For a free algebra $H,$ the module ${\mathcal{M}}\subset V$ has simple form
in the (minimal) Hirsch resolution $(RH,d),$ i.e., ${\mathcal{M}}^{<0,\ast
}=0.$ Indeed, given an odd dimensional multiplicative generator $x\in H$ and
a representative $x_{0}\in R^{0}H$ of $x,$ the elements $x_{n}$ in the
sequence (\ref{sms}) can be defined as $x_{n}=\frac{(-1)^{n}}{(n+1)!}
\,x_{0}^{\smallsmile _{1}n+1}$ and hence $x_{n}\in {\mathcal{E}}$ for $n\geq
1.$ In particular, there is a map of dg algebras $(RH,d)\rightarrow A$ and
hence an isomorphism of dg coalgebras $H^{\ast }(BA)\approx H^{\ast }(BH)$
for a dga $A$ with $H=H^{\ast }(A)$ (a free $\Bbbk$-algebra $H$ is intrinsically $
\Bbbk $-formal). Regarding the filtered Hirsch model $(RH,d_{h}),$ the
perturbation $h$ may be non-zero only on $\mathcal{T}.$ More precisely,
according to Example \ref{hoh} the cohomology class $[h(a\cup _{2}b)]\in
H^{\ast }(RH,d_{h})$ is defined by $\rho a\ast \rho b\in H$ for $a,b\in
V^{0,\ast }.$ Since $H^{\ast }(BH)\approx H^{\ast }(BA)\approx H^{\ast }(
\bar{V},\bar{d}_{h})$ (cf. (\ref{Vbar})), the multiplication $\mu _{E}^{\ast
}$ on $H^{\ast }(BH)$ is induced by the $\smile _{1}$-product on $V$ (cf.
Remark \ref{remark3}). Therefore, the Lie bracket on $H^{\ast }(BH)$ is
determined by the bracket
\begin{equation*}
\lbrack a,b]=a\smile _{1}b-(-1)^{(|a|+1)(|b|+1)}b\smile _{1}a
\end{equation*}%
on $V.$ The observation that $s^{-1}[a,b]$ is cohomologous to $s^{-1}h(a\cup
_{2}b)$ in $\bar{V}$ for all $a,b\in V^{0,\ast }$ completes the proof.
\end{proof}

\begin{remark}
\label{lie} Note that the transgressive component $h^{tr}$ evaluated on the
elements $a_{1}\cup _{2}\cdots \cup _{2}a_{n}\in \mathcal{T}$ for $a_{i}\in
V^{0,\ast },\,n\geq 3,$ determines higher order operations on $Tor^{A}(\Bbbk
;\Bbbk )$ that extend the Lie algebra structure to an $L_{\infty }$-algebra
structure.
\end{remark}

For example, a polynomial algebra $P=\Bbbk \lbrack x_{1},...,x_{n}]$
provides the case of $H^{\ast }$ in the theorem. Indeed, in general, to
calculate the Hochschild cohomology of an algebra $P$ construct a small
complex $(C_{V}^{\bullet }(P),\bar{d}),$ which is quasi-isomorphic to $
C^{\bullet }(P;P)$ as follows (compare \cite{J-M.hoh}): Fix an ordinary
multiplicative resolution $\rho :RP\rightarrow P$ with $RP=T(V),$ view $P$
as an $RP$-bimodule via $\rho ,$ and let $B(\rho )^{\bullet }:C^{\bullet
}(P;P)\rightarrow C^{\bullet }(RP;P)$ be a quasi-isomorphism induced by $
B(\rho ):B(RP)\rightarrow BP.$ Set $(C_{V}^{\bullet }(P),\bar{d})=(Hom(\bar{V
},P),\bar{d})$ in which $\bar{d}$ is defined for $f\in C_{V}^{\bullet }(P)$
by $\bar{d}f=g,$
\begin{equation*}
g(\bar{x})=\underset{1\leq i\leq k}{\sum }(-1)^{\nu _{i}}\rho (v_{1})\cdots
f(\bar{v}_{i})\cdots \rho (v_{k}),\,\,\,\,dx=\sum v_{1}\cdots
v_{k},\,v_{i}\in V,\,k\geq 1,
\end{equation*}
$\nu _{i}=(|f|+1)(|v_{1}|+\cdots +|v_{i-1}|),$ and define a chain map $\chi
:C_{V}^{\bullet }(P)\rightarrow C^{\bullet }(RP;P)$ by $\chi f=f^{\prime },$
\begin{equation*}
f^{\prime }(\bar{x})=\left\{
\begin{array}{lll}
f(\bar{x}), & x\in V,\vspace{1mm} &  \\
\underset{1\leq i\leq n}{\sum }(-1)^{\nu _{i}}\rho (v_{1})\cdots f(\bar{v}
_{i})\cdots \rho (v_{n}), & x=\sum v_{1}\cdots v_{n},\,v_{i}\in V,n\geq 2. &
\end{array}
\right.
\end{equation*}
Isomorphism (\ref{Vbar}) implies that $\chi $ is a homology isomorphism. On
the other hand, the $\smile $-product on $C^{\bullet }(P;P)$ induces a $
\smile $-product on $C_{V}^{\bullet }(P);$  more precisely,
we have that $\bar V$ is a coalgebra with
 the coproduct
$\bar{\Delta}:\bar{V}\rightarrow \bar{V}\otimes \bar{V}$ induced by the
standard coproduct of $BP$
and, consequently, $Hom(\bar V,P)$ is endowed with the standard  $\smile$-product.
 When $P$ is polynomial, the minimal $V^{\ast }$
can be thought of as generated by the iterations of a (commutative) $\smile
_{1}$-product (\cite{saneCLASS}); consequently, $(\bar{V}^{\ast },\bar{\Delta
})$ is an exterior coalgebra. Dually, $\bar{V}_{\ast }$ is an exterior
algebra on generators $\bar{x}_{1},...,\bar{x}_{n}.$ Furthermore, $\bar{d}=0$
and hence $H(C_{V}^{\bullet }(P),\bar{d})=C_{V}^{\bullet }(P).$ Thus the
Hochschild cohomology $H^{\ast }$ is isomorphic to the algebra $
C_{V}^{\bullet }(P)\approx \bar{V}_{\ast -1}\otimes P^{\ast },$ which is the
tensor product of an exterior algebra and a polynomial algebra, as required.

\subsection{Symmetric Massey products in $C^*(X;\Bbbk)$ and powers in the
loop homology $H_*(\Omega X;\Bbbk)$}

Let $A_{\ast }$ be a dg coalgebra over a field $\Bbbk $ and let $A^{\ast
}=Hom(A_{\ast },\Bbbk )$ be a dg algebra so that $H(A^{\ast })=Hom(H(A_{\ast
}),\Bbbk ).$ Let
\begin{equation*}
\iota :H(BA^{\ast })\rightarrow Hom(H(\Omega A_{\ast }),\Bbbk )),
\end{equation*}
be the canonical map, where $\Omega A_{\ast }$ denotes the cobar construction of the coalgebra $
A_{\ast }.$ Given the suspension map $\sigma ^{\ast }:H^{\ast }(A^{\ast
})\rightarrow H^{\ast -1}(BA^{\ast }),$ let $x\in H_{\ast }(A^{\ast })$ and $
y\in H_{\ast -1}(\Omega A_{\ast }),$ where $y$ is a basis element with $
\iota (\sigma ^{\ast }x)(y)=1\in \Bbbk ,$ and $\iota (\sigma ^{\ast
}x)(y^{\prime })=0$ for any basis element $y^{\prime }\neq y.$

Suppose that $\langle x\rangle ^{n}$ is defined for $x.$ Let $
\{a_{i}\}_{0\leq i<n}$ be a defining system of $\langle x\rangle ^{n}$ with $
a_{0}\in A^{\ast }$ a representative cocycle of $x.$ Then $\bar{a}_{0}\in
BA^{\ast }$ is a cocycle with $[\bar{a}_{0}]=\sigma ^{\ast }x$ and $
\{a_{i}\}_{0\leq i<n}$ lifts to a cocycle $a\in BA^{\ast }$ so that the
cohomology class $[a]\in H^{\ast }(BA^{\ast })$ is represented by the $y^{n}$
(the $n^{th}$-power of $y)$ in $H_{\ast }(\Omega A_{\ast })$ via the map $
\iota .$ Then Theorem \ref{mzero} immediately implies the following:

\begin{theorem}
Let $X$ be a simply connected space, let $\Bbbk $ be a field of
characteristic zero, and let $\sigma _{\ast }:H_{\ast }(\Omega X;\Bbbk
)\rightarrow H_{\ast +1}(X;\Bbbk )$ be the suspension map. If $y\in H_{\ast
}(\Omega X;\Bbbk )$ such that $y\notin \operatorname{Ker}\sigma _{\ast }$ and $%
y^{2}\neq 0$, then $y^{n}\neq 0$ in $H_{\ast }(\Omega X;\Bbbk )$ for all $%
n\geq 2.$
\end{theorem}

Finally,
recalling the  connection between symmetric Massey products
and twisting elements in $A^{\ast },$ which arise from the sequences $
\{a_{i}\}_{i\geq 0}$ above, we remark that the observation above relates the
existence of twisting elements in $A^{\ast }$ with the existence of
polynomial generators in $H_{\ast }(\Omega A_{\ast }).$

\vspace{5mm}

\end{document}